\documentclass[singlecolumn]{elsart}    

\usepackage{algorithm}
\usepackage{algpseudocode}
\usepackage[table]{xcolor}
\usepackage{graphicx}
\usepackage{epsfig}
\usepackage{amssymb}
\usepackage{amsfonts}

\usepackage{epstopdf}
\usepackage{bbm}
\usepackage{multirow}
\usepackage{tikz}
\usetikzlibrary{patterns, shapes}
\usepackage{wrapfig}
\usepackage{amsmath}

\usepackage{caption}
\DeclareOldFontCommand{\cal}{\normalfont\rmfamily}{\mathcal}
\newcommand{\nn}{\nonumber}

\newcommand{\mz}{\color{black}}
\newcommand{\tr}{\operatorname{tr}}

\newcommand{\Rs}{\mathbb R}
\newcommand{\Ds}{\mathcal D}
\newcommand{\Ls}{\mathcal L}

\newcommand{\Bs}{\mathcal B}
\newcommand{\Es}{\mathbb{E}}

\newcommand{\al}[1]{\begin{align} #1 \end{align}}

\begin{document}

\begin{frontmatter}

\title{An update-resilient Kalman filtering approach} 

\thanks[footnoteinfo]{This paper was not presented at any IFAC
meeting. Corresponding author: S. Yi.}

\author[Italy]{Shenglun Yi}\ead{shenglun@dei.unipd.it},    
\author[Italy]{Mattia Zorzi}\ead{zorzimat@dei.unipd.it},               

\address[Italy]{Department of Information Engineering, University of Padova, Via Gradenigo 6/B, 35131 Padova, Italy}  

\begin{keyword}                           
Robust Kalman filtering; Kullback-Leibler divergence; Minimax game.         
\end{keyword}                             

\begin{abstract}                          
We propose a new robust filtering paradigm considering the situation in which model uncertainty, described through an ambiguity set, is present only in the observations.   We derive the corresponding robust estimator,  referred to as update-resilient Kalman filter, which appears to be novel compared to existing minimax game-based filtering approaches.  Moreover, we characterize the corresponding least favorable model and analyze the filter stability. Finally, some numerical examples show the effectiveness of the proposed estimator.
\end{abstract}

\end{frontmatter}

 \section{Introduction}
Various  state estimation  paradigms have been proposed in the literature  to tackle model  uncertainty such as considering parametric approaches \cite{ROCHA2021105034,GHAOUI_CALAFIORE_2001,9113473},  filters robust to outliers \cite{huang2016robust,huang2019robust} and minimax approaches \cite{YOON_2004,MOP_23,ROBUSTNESS_HANSENSARGENT_2008,ROBUST_STATE_SPACE_LEVY_NIKOUKHAH_2013}.
In the latter, uncertainty is characterized by an ambiguity set that captures the ``mismatch'' between the actual and nominal  models. In particular, to allow the uncertainty to be uniformly distributed over time, a well-established paradigm is to define this set  at each time step which is a ball whose radius represents the level of uncertainty  and its center corresponds to the nominal model \cite{ROBUST_STATE_SPACE_LEVY_NIKOUKHAH_2013}.
The minimax paradigm is formulated as a dynamic game: the maximizer selects the least favorable model within the prescribed ambiguity set, while the other player, i.e. the estimator, seeks to minimize the estimation error with respect to the chosen maximizer. Various extensions of this paradigm have been proposed, including ambiguity  sets defined using different metrics \cite{10654520,STATETAU_2017}, formulations with degenerate densities \cite{yi2020low,yi2021robust},  adaptive formulations \cite{yi2024robust} and   the case of nonlinear state space models   \cite{longhini2021learning}.  Finally, the relaxation of such minimax games lead to the so called risk sensitive filters \cite{H_INF_HASSIBI_SAYED_KAILATH_1999,speyer_2008,huang2018distributed}.

These minimax game-based filters share a common feature: the uncertainties are captured by an ambiguity set defined over the  entire state space model, meaning that uncertainties affect both the state dynamics and the observations.  This particular feature  has a significant impact on the structure of the resulting robust estimators.  The latter exhibit a structure similar to that of the Kalman filter  (KF), in which the ``resilience'' to uncertainty is in the prediction stage, that is, the covariance matrix of the prediction error is modified.
An exception, however, is represented by \cite{abadeh2018wasserstein}, where the ``resilience'' is applied in the update stage.  The minimax formulation of this robust ﬁltering problem is based on an ambiguity set defined through the Wasserstein distance. Nevertheless, that paper does not provide an  explicit  interpretation of the underlying uncertainty in terms of  the state space model. Furthermore,  it lacks the characterization of the state space representation corresponding to the least favorable scenario, and  does not guarantee the filter stability even when the tolerance (i.e. the Wasserstein radius) is small.

In many real-world systems,  the state dynamics are relatively well-understood and accurately modeled, while  the uncertainty is mainly in the observations. An example is represented by the displacement estimation problem of a mass spring damper system, where the displacement is measured by sensors that are susceptible to various  uncertainties. In such unbalanced uncertainty scenarios, where the ``dirty'' sensor data constitutes the dominant source of uncertainty, it becomes important to reconsider whether characterizing the ambiguity set in terms of the entire state space model remains the most appropriate choice. Indeed, this characterization may still appear reasonable from an intuitive standpoint,  since in practice there exists a mismatch (even though mild) between the actual and nominal state process models. Thus, an important question concerns whether the robust filters like \cite{ROBUST_STATE_SPACE_LEVY_NIKOUKHAH_2013,10654520}, which are resilient in the prediction stage, can perform effectively in the presence of unbalanced uncertainty.

%

%

In this paper, we propose a new robust state space filtering problem in which uncertainties arise in the observation model.  To this end, we define a novel ambiguity set that pertains solely to the density characterizing the observation model  by means of the Kullback-Leibler (KL) divergence. Our analysis shows that the resulting estimator exhibits a Kalman filter-like structure, in which the resilience is  in the update stage.  
  This paradigm is fundamentally different from the ones in \cite{ROBUST_STATE_SPACE_LEVY_NIKOUKHAH_2013,10654520}.   {\mz  This fact is also reflected in the corresponding least favorable model, whose approximation through a state-space representation exhibits a larger dimension than those in }\cite{ROBUST_STATE_SPACE_LEVY_NIKOUKHAH_2013,10654520}.    Moreover, we analyze the filter stability in the worst case scenario.  The numerical results show that even in the presence of a mild mismatch between the actual and nominal state process models, the proposed estimator still outperforms those in \cite{ROBUST_STATE_SPACE_LEVY_NIKOUKHAH_2013,10654520}. This is because in such unbalanced uncertainty scenarios, characterizing the ambiguity set in terms of the entire state space  model result in overly conservative estimates, as modelers has opportunity to allocate the ``mismatch'' budget uniformly across both the state process and the observation model.
Furthermore, numerical results show that the proposed estimator performs similarly to the one in \cite{abadeh2018wasserstein}. Indeed, since both are resilient in the update stage,   our  theoretical framework  suggests that the estimator proposed in \cite{abadeh2018wasserstein} implicitly postulates uncertainty only in the observation model.  In addition, the experimental evidence  indicates that our approach is significantly more computationally efficient than \cite{abadeh2018wasserstein}, i.e. a property that is particularly important for real-time applications.
Finally, we consider the relaxed version of the aforementioned  minimax game and show that the resulting estimator is a new risk sensitive filter.

The outline of the paper is as follows. In Section \ref{sec_2} we introduce the problem formulation.  In Section \ref{sec_3} we derive the update-resilient Kalman filter and its corresponding least favorable model. In Section \ref{sec_4} we analyze
the stability of the proposed filter. In Section \ref{sec_ne} we provide some numerical examples.
In Section \ref{sec_rsf} we derive the corresponding risk sensitive estimation paradigm. Finally, in Section \ref{sec_con} we draw the conclusions.



\emph{Notation}: Given a symmetric matrix $K: K > 0~ (K \geq 0)$ means that $K$ is positive
(semi-)deﬁnite;  $\sigma_{max}(K)$ is the maximum eigenvalue of $K$;  $\tr(K)$ and $\det(K)$ denote the trace and determinant of $K$, respectively. Then, $v\sim \mathcal N(m,R)$ means that $x$ is a Gaussian random vector with mean $m$ and covariance matrix $R$.   Finally, throughout the paper,
we omit the arguments of probability  density functions, e.g. $\psi_t$ instead of $\psi_t(y_t|x_t)$,  when the context is clear in order to  keep the notation concise.  Moreover, $\hat{x}_{t|t}$ denotes the posterior state estimate, while $\hat{x}_{t}$ denotes the prior state estimate. The corresponding covariance matrices follow the same subscript notation: the symbol used for the covariance with subscript $t$ refers to the prior covariance matrix, and the same symbol with subscript $t|t$ refers to the posterior covariance matrix.


\section{Problem formulation} \label{sec_2}
We consider a robust  discrete-time state space filtering problem where perturbations occur solely in the observations. In this scenario, the actual model for the state process  is assumed to be known and takes the form:
\begin{equation}\label{nomi_linea_mod_proc}
x_{t+1}= A x_t + \epsilon_t
 \end{equation}
where $ A \in \mathbb{R}^{n\times n}$,  $x_t \in \mathbb{R}^{n}$ is the state vector  and
$\epsilon_t \in \mathbb{R}^{n}$ is white Gaussian  noise  with the zero mean and covariance matrix $ Q \in \mathbb{R}^{n\times n}$.  Thus,  model (\ref{nomi_linea_mod_proc}) can be entirely  characterized
by the following transition density ${p_{t}}(x_{t+1}|x_t) \sim \mathcal{N} (A x_t , ~ Q)$ and the initial state  $x_0 \sim p_0(x_0)$, which is assumed to be independent  from $\epsilon_t$. The nominal model for the observations is
\begin{equation}\label{nomi_linea_mod_obv}
y_{t}=  C x_t + \varepsilon_t
 \end{equation}
where $ C \in \mathbb{R}^{m\times n}$,  $y_t \in \mathbb{R}^{m}$ is the observation vector, and $\varepsilon_t \in \mathbb{R}^{m}$ is white  Gaussian  noise  with zero mean and   covariance matrix $ R \in \mathbb{R}^{m\times m},$ which is independent from $\epsilon_t$. We also assume that   $ Q >0$ and $R >0$. The nominal  model (\ref{nomi_linea_mod_obv}) is completely  characterized
by  the nominal density $ \psi_t(y_t|x_t) \sim \mathcal{N} ( C x_t , ~ R).$  Accordingly,  the state space model (\ref{nomi_linea_mod_proc})-(\ref{nomi_linea_mod_obv}) over the finite  time interval  $t\in\{0\ldots N\}$ is characterized  by the following nominal joint probability density of  $  X_N :=\{{x_{0}} {\ldots}  {x_{N+1}}\}$ and $Y_N :=\{{y_{0}}  {\ldots}  {y_{N}}\}$:
 \begin{equation*}
p\left(X_N,Y_{N}\right)=p_0\left(x_{0}\right) \prod_{t=0}^{N} {p_{t}}\left(x_{t+1} | x_{t}\right) \psi_{t}\left(y_{t} | x_{t}\right).
\end{equation*}
Let $\tilde \psi_{t}\left(y_{t} | x_{t}\right)$ denote the conditional density  characterizing the actual model for the observations.  Note that,  $\tilde \psi_t$ is not required to follow any specific form.
We assume that  actual  density of  $X_{N},Y_N$, say  $\tilde{p}(X_N,{Y}_{N})$, follows a Markov structure similar to the nominal one:
\begin{equation}\label{tildepZ}
\tilde{p} \left(X_N,Y_{N}\right)={p}_0\left(x_{0}\right) \prod_{t=0}^{N} {p_{t}}\left(x_{t+1} | x_{t}\right)  \tilde \psi_{t}\left(y_{t} | x_{t}\right).
\end{equation}

We measure the modeling mismatch between the nominal density of  $X_N, Y_N$ and the actual one through the KL divergence  \cite{kullback1951information}:
\begin{equation}\label{KL_p}
\Ds_{KL}(\tilde{p}, p)=\int \tilde{p}  \ln \left(\frac{\tilde{p}}{p}\right) d X_N dY_N.
\end{equation}
Taking the expectation of
\begin{equation*}
\ln \left(\frac{\tilde{p}}{p}\right)=\sum_{t=0}^N \ln \left(\frac{\tilde \psi_t\left(y_t \mid x_t\right)}{\psi_t\left(y_t \mid x_t\right)}\right)
\end{equation*}
with respect to  $\tilde p({X}_{N},Y_N),$ we see the KL divergence (\ref{KL_p}) takes the form:
\begin{equation}\nn \begin{aligned}
  \Ds_{KL}(\tilde{p},  p)
 = \sum_{t=0}^{N} \Ds_{KL}(\tilde \psi_t,\psi_t)\end{aligned}
\end{equation}
where
\begin{align}\label{DS_O}
\Ds_{KL}(\tilde \psi_t,\psi_t)&:= \tilde \Es\left[\ln\left(\frac{\tilde \psi_t}{\psi_t}\right)\right]\nn\\
& =\iint  \tilde \psi_t(y_t | x_t)  p_t(x_t) \ln\left(\frac{\tilde \psi_t}{\psi_t}\right)d y_t d x_t,
\end{align}
and $p_t(x_t)$ denotes the marginal density of $x_t$. At time $t$, given the observations $Y_{t-1}:=\{y_s, ~ s\leq t-1\}$, we  assume that the   actual  density $\tilde \psi_t(y_t|x_t)$  belongs to the following convex ambiguity set\footnote{The adjective ambiguity is used to express the lack of the precise knowledge about the actual model.}:
\begin{equation}\label{ballpsi}
{\mathcal B}_t:=\left\{\, \tilde{\psi}_{t} \hbox{ s.t. } \Ds(\tilde \psi_t,\psi_t) \leq c_{t}\right\}
\end{equation}
which can be regarded as a ``ball'', with tolerance $c_t>0$ representing its radius about the nominal  density $\psi_t$. 
Such ball is with respect to the metric induced by the KL divergence (\ref{DS_O}) conditioned on the observations  up to time $t$, i.e.
\begin{align*}
 \Ds(\tilde \psi_t,\psi_t) &:=  \tilde \Es \left[\ln\left(\frac{\tilde \psi_t}{\psi_t}\right)\Bigg| Y_{t-1}\right]\nn\\
& =\iint  \tilde \psi_t(y_t | x_t) \tilde p_t(x_t|Y_{t-1}) \ln\left(\frac{\tilde \psi_t}{\psi_t}\right)d y_t d x_t
\end{align*}
where $\tilde p_t\left({x}_t | Y_{t-1}\right)$  is  the actual \emph{a priori}
conditional density  of $x_t$ conditioned on $Y_{t-1}$.  The latter is, in general, different from the nominal one, say $ p_t\left({x}_t | Y_{t-1}\right)$, because the conditioning depends on the model for the observations which is  affected by uncertainty.  In practice, the tolerance is unknown. One possible strategy is to learn a constant tolerance  from data, i.e. $c_t=c$. More precisely, given a dataset of observations and a finite set of candidate values for $c$, the tolerance is estimated as the candidate corresponding to the robust output predictor that minimizes the sample variance of the output prediction error, \cite{11107768}.
Notice that, the ambiguity set \eqref{ballpsi} differs from the one considered in \cite{ROBUST_STATE_SPACE_LEVY_NIKOUKHAH_2013}. While we consider uncertainty only in the measurement equation \eqref{nomi_linea_mod_obv}, \cite{ROBUST_STATE_SPACE_LEVY_NIKOUKHAH_2013} considers uncertainty in both the state process and measurement equations.

The unbalanced uncertainty scenario described above arises in many practical problems. Indeed, the state process model is often assumed to be known, since it is typically derived from well-established physical laws. By contrast, although uncertainty in the observation model can be reduced through the use of advanced sensors, such devices are often very expensive. As a result, less advanced sensors are typically employed in practice due to budgetary constraints. Moreover, sensor uncertainty cannot be entirely eliminated, even with advanced hardware, because of unexpected faults, such as mechanical stress or calibration drift, or environmental factors, such as multipath effects in GPS and reflections in LiDAR.

{\textbf{The significance of the ambiguity set.}} Consider (\ref{nomi_linea_mod_obv}) as a regression model where $x_t$  acts as the input of such regression model, whose prior is defined by (\ref{nomi_linea_mod_proc}), and $y_t$ is the output.
Note that,  the latter is completely described by $\psi_t$.
Assume to collect a set of $M$ independent  state-output data corresponding to time $t$, i.e. $\mathbf D_{t,M}:=\{ (y_t^k,x_t^k), \; k=1\ldots M \}$. We assume these data are generated by the actual model with $Y_{t-1}=\{y_{0}\ldots y_{t-1}\}$ fixed.
The log-likelihood of the data based on the regression model  (\ref{nomi_linea_mod_obv}) is
\al{\ell(\mathbf D_{t,M}; \psi_t):=\sum_{k=1}^M \ln \psi_t (y_t^k|x_t^k).\nn} Thus, the expected log-likelihood
is
\al{\ell_\infty(\psi_t)&:=\lim_{M\rightarrow \infty} \frac{1}{M}\ell(\mathbf D_{t,M}; \psi_t)= \tilde \Es \left[\ln \psi_t | Y_{t-1}\right]\nn\\
& =\iint  \ln\left({\psi_t}\right)\tilde \psi_t(y_t | x_t) \tilde p_t(x_t|Y_{t-1}) d y_t d x_t \nn\\
& =  - \iint  \tilde \psi_t(y_t | x_t) \tilde p_t(x_t|Y_{t-1}) \ln\left(\frac{\tilde \psi_t}{\psi_t}\right)d y_t d x_t \nn\\ 
&  ~~~~ + \iint  \ln({\tilde \psi_t})\tilde \psi_t(y_t | x_t) \tilde p_t(x_t|Y_{t-1}) d y_t d x_t \nn \\
&=\kappa_t- \Ds(\tilde \psi_t,\psi_t) \nn}
where the limit above almost surely exists (i.e. equality holds with probability one) and
\al{\kappa_t:=\iint  \ln({\tilde \psi_t})\tilde \psi_t(y_t | x_t) \tilde p_t(x_t|Y_{t-1}) d y_t d x_t \nn}
is a term not depending on the regression model (\ref{nomi_linea_mod_obv}).

We conclude that the ambiguity set (\ref{ballpsi})  includes models whose expected log-likelihood is bounded below by a threshold determined by the tolerance $c_t$: $$\ell_\infty(\tilde \psi_t)\geq \kappa_t-c_t.$$
It is worth noting that the ambiguity set in (\ref{ballpsi}) is fundamentally different from the one proposed in \cite{ROBUST_STATE_SPACE_LEVY_NIKOUKHAH_2013,10654520}, i.e. {\small\al{\label{oldAB}\check{\mathcal B}_t:=\left\{(\psi_t,\tilde p_t) \hbox{ s.t. } \int \int \int \tilde \psi_t\tilde p_t \ln \left(\frac{\tilde \psi_t\tilde p_t}{\psi_tp_t}\right)d y_t d x_t dx_{t+1}\right\},}}which considers uncertainty in both $\psi_t(y_t|x_t)$ and $p_t(x_{t+1}|x_t)$. However, the latter is an inappropriate choice when the primary source of uncertainty lies in $\psi_t$, i.e. the measurement equation. In such cases, designing a minimax estimator that allows uncertainty in both $\psi_t(y_t|x_t)$ and $p_t(x_{t+1}|x_t)$ grants the hostile player the freedom to allocate  a relevant portion of the  mismatch budget to the state dynamics (\ref{nomi_linea_mod_proc}), which misrepresents the actual scenario. As a result, the performance of the estimator will be severely compromised.

The aim of this paper is to address the following problem.
\begin{prob}
Design a robust estimator of $x_t$ given 
$Y_t$ that minimizes the upper bound on the variance of the filtering error for all models whose conditional density of the observation model belongs to the ambiguity set (\ref{ballpsi}).
\end{prob}

In what follows, to ease the exposition, we will consider the case in which the state space
model (\ref{nomi_linea_mod_proc})-(\ref{nomi_linea_mod_obv}) is time-invariant, however the results we will present can be straightforwardly extended to time-varying case.

\section{Update--resilient Kalman filter}\label{sec_3}
We design the robust estimator of $x_{t}$ given $Y_t$ as the solution to the following dynamic minimax game:
\begin{equation}
 \label{minimax_psi} (\tilde \psi^\star_t, g^\star_t) =\arg \underset{g_t \in \mathcal{G}_{t}}{\mathrm{min}}\max_{\tilde{\psi}_{t} \in {\mathcal{B}}_{t}} J_t(\tilde {\psi}_t,g_t)
\end{equation}
where the objective function is given by  \begin{equation}
 \begin{aligned}
J_t&(\tilde {\psi}_t,g_t)=\frac{1}{2}\tilde{\mathbb{E}}\left[\left\|   x_{t}-g_t\left(y_{t}\right) \right\|^{2} | {Y}_{t-1}\right]\\
&=\frac{1}{2} \iint \left\|  x_{t}-g_{t}\left(y_{t}\right) \right\|^{2} \tilde{\psi}_{t}(y_t | x_t)   {\tilde p}_{t}\left(x_{t} | Y_{t-1}\right) d y_{t} d x_{t}\nn
\end{aligned}
\end{equation} which represents the variance of the filtering error under the actual density $\tilde \psi_t$. The latter is restricted to lie within the ambiguity set $\mathcal{B}_{t}$, i.e. we impose an upper bound on the mismatch between $\tilde{\psi}_t$ and the nominal density $\psi_t$ which is measured using the conditional KL divergence $\mathcal{D}(\tilde{\psi}_t, \psi_t)$. $\mathcal{G}_{t}$ denotes the class of estimators with finite second-order moments with respect to all the densities $\tilde{\psi}_{t} \tilde{p}_t(x_t|Y_{t-1})$ such that $\tilde{\psi}_{t} \in {\mathcal{B}}_{t}$.   Note that, $\tilde \psi_t$ must satisfy also the normalization constraint:
\begin{equation}
I_{t} (\tilde{\psi}_{t} )  := \iint  \tilde{\psi}_{t}(y_t | x_t) {\tilde p}_{t}\left(x_{t} | Y_{t-1}\right)  d y_{t} d x_{t}=1. \label{I}
\end{equation}
Note that, $J_t$ is quadratic in $g_t$, in particular it is convex in $g_t$, and linear in $\tilde \psi_t$. Then, on the basis of the Von Neumann’s minimax theorem  \cite{aubin2006applied}, there exists a saddle point $(\tilde \psi^\star_t,~g^\star_t)$ such that:
\begin{equation}\label{J_minimax2}
                    J_t(\tilde \psi_t,~ g^\star_t) \leq J_t(\tilde \psi^\star_t,~g^\star_t) \leq J_t(\tilde \psi^\star_t,~ g_t)
                  \end{equation}
for any $g_t \in \mathcal{G}_{t}$ and $\tilde \psi_t \in \mathcal{B}_{t}.$ From the second inequality of (\ref{J_minimax2}), it follows that the minimizer of (\ref{minimax_psi}) coincides with the expectation of $x_t$ taken with respect to the actual filtering density, i.e. $
\tilde{p}_{t}\left(x_{t} | Y_{t}\right),
$ which depends on  $\tilde \psi^\star_t$. Next,  we characterize the  solution to the  minimax problem (\ref{minimax_psi}).

\begin{lem}
\label{lemma1}
Under the assumption that $\tilde p(x_t|Y_{t-1})$ is different from zero almost everywhere,  the maximizer of (\ref{minimax_psi})
takes the form:
\begin{equation} \label{psi_0}
\tilde{\psi}^0_{t} (y_t|x_t) =\frac{1}{M_t (\lambda_{t}) } \exp \left( \frac{1}{2\lambda_t}  \left\|  x_{t}-g_t(y_t) \right\|^{2}\right) \psi_{t} (y_t|x_t)
\end{equation}
where   \begin{equation}\label{Mt}
\scalebox{0.9}{$
\begin{aligned} M_t(\lambda_{t})=\iint \hspace{-0.1cm} \exp \left( \hspace{-0.05cm}\frac{1}{2\lambda_{t}}  \left\|  x_{t}-g_t(y_t) \right\|^{2}\hspace{-0.05cm} \right)\hspace{-0.1cm}  \psi_{t} \tilde p_{t}\left(x_{t} | Y_{t-1}\right)   d y_{t} d x_{t}  \end{aligned} $}
\end{equation} is the normalizing constant such that (\ref{I}) holds. Moreover,  $\lambda_{t}>0$ is the unique solution to the equation   $\Ds(\tilde \psi^0_{t}, \psi_{t})=c_t $.
\end{lem}
\begin{pf}
%
We want to maximize $J_{t} (\tilde{\psi}_{t} , g_t) $  with respect to $\tilde \psi_t \in \Bs_t$ using the Lagrangian multipliers theory.  More precisely, we consider the  Lagrange function:
\begin{equation}\label{lag}
\begin{aligned}
\Ls &(\tilde{\psi}_{t}, ~ \lambda_{t},~ \beta_t) \\
&=  J_{t} (\tilde{\psi}_{t}, g_t) + \lambda_{t} ( c_{t} - \Ds(\tilde{\psi}_{t}, {\psi}_{t}))  +\beta_t (I_{t}(\tilde{\psi}_{t})-1) \\
&=  \iint \left(  \frac{1}{2} \left\|  x_{t}-g_{t}\left(y_{t}\right) \right\|^{2} - \lambda_{t} \ln\left(\frac{\tilde \psi_t}{\psi_t}\right) + \beta_t \right) \\
 &\hspace{0.5cm} \times \tilde{\psi}_{t}   {\tilde p}_{t}\left(x_{t} | Y_{t-1}\right) d y_{t} d x_{t} + \lambda_{t} c_t - \beta_t
\end{aligned}
\end{equation}
where $\lambda_{t} > 0$ is the  Lagrange multiplier corresponding to constraint $\tilde \psi_t \in \Bs_t$ and
$\beta_t$ is the  one corresponding to constraint (\ref{I}). Notice that, the Lagrangian (\ref{lag}) is concave with respect to $\tilde \psi_t$  because it is the sum of two linear functions in $\tilde \psi_t$, i.e. $J_{t}$ and $I_{t}$, and $-\Ds(\tilde{\psi}_{t}, {\psi}_{t})$ which is concave with respect to the second argument, \cite{COVER_THOMAS}. We show that  $\Ls_{t}$   has a unique stationary point,  then the latter is the unique point of maximum for $\Ls_{t}$. The first variation of $\Ls_t$  along the direction $\delta \tilde \psi_t$ is given by
\begin{equation*}
\begin{aligned} \delta & \Ls(\tilde{\psi}_{t}, \lambda_{t}, \beta_{t} ;\delta \tilde{\psi}_{t})
= \iint \delta \tilde{\psi}_{t}\tilde{p}_{t}\left(x_{t} | Y_{t-1}\right) \\
& \times( \frac{1}{2} \left\|  x_{t}-g_t(y_t)\right\|^{2} +  \lambda_{t} \ln\left(\frac{{\psi}_{t}}{\tilde{\psi}_{t}}\right) + \beta_t -\lambda_{t})   d y_t d x_t.
\end{aligned}
\end{equation*}
Accordingly, the stationary point $\tilde \psi^0_t$ must satisfy $$\delta \Ls(\tilde{\psi}^0_{t}, \lambda_{t}, \beta_{t} ;\delta \tilde{\psi}_{t}) = 0 $$ for any function $\delta \tilde \psi_t$.  Since $\tilde p(x_t|Y_{t-1}) $ is  different from zero almost everywhere, we obtain
$$ \frac{1}{2} \left\|   x_{t}-g_t(y_t)\right\|^{2} +  \lambda_{t} \ln\left(\frac{{\psi}_{t}}{\tilde{\psi}^0_{t}}\right) + \beta_t - \lambda_{t} = 0.$$
Then, it is easy to see that
$$\tilde{\psi}^0_{t}(y_t|x_t) = {\psi}_{t} \exp \left(\frac{1}{2 \lambda_{t}}  \left\|   x_{t} -g_t(y_t)\right\|^{2} +  \frac{ \beta_t}{\lambda_t} -1 \right).$$
It is not difficult to see that the optimal value for $\beta_t$, say $\beta_t^0$, is such that
$\exp( \beta^0_t/\lambda_{t} -1) = M_t(\lambda_{t})^{-1}$
where  $M_t(\lambda_{t})$ is defined in (\ref{Mt}). Thus, we obtain (\ref{psi_0}).

It remains to consider the dual problem for the  Lagrange multiplier  $\lambda_t  > 0$. More precisely,  the dual function is given by:
\begin{equation} \nn
\begin{aligned}
\tilde{\Ls}(\lambda_t) &= \Ls (\tilde{\psi}^0_{t},  \lambda_t, \beta^0_t) \\
&=  \iint \left(  \frac{1}{2} \left\|  x_{t}-g_{t}\left(y_{t}\right) \right\|^{2} - \lambda_t \ln\left(\frac{\tilde \psi^0_t}{\psi_t}\right)  \right) \\
 &\hspace{2cm} \times \tilde{\psi}_{t}   {\tilde p}_{t}\left(x_{t} | Y_{t-1}\right) d y_{t} d x_{t} + \lambda_t c_t \\
& =  \lambda_t(\ln(M_t(\lambda_t)) + c_t) .
\end{aligned}
\end{equation}
Notice that,  $ \exp\left(\frac{1}{2\lambda_t} \left\|  x_{t}-g_t(y_t) \right\|^{2}\right) \rightarrow  1$
as $\lambda_t \rightarrow \infty$.
Thus, $\ln(M_t(\lambda_t)) \rightarrow  0,$ so $\tilde{\Ls}(\infty) = \infty $ since $c_t>0$.
 Accordingly, the infimum cannot be attained for  $\lambda_t \rightarrow \infty$. Finally,
we consider:
$$\begin{aligned}
\frac{\mathrm{d}  }{\mathrm{d} \lambda_t} \tilde{\Ls} (\lambda_t) &= \ln(M_t(\lambda_t)) + \lambda_t M^{-1}_t(\lambda_t) \frac{\mathrm{d}  }{\mathrm{d} \lambda_t} M_t(\lambda_t) +c_t\\
&=c_t - \lambda^{-1}_t J_t(\tilde {\psi}^0_t(\lambda_t),g_t) + \ln(M_t(\lambda_t)) \\
&=c_t - \Ds(\tilde{\psi}^0_{t}(\lambda_t), {\psi}_{t}).
\end{aligned}
 $$
In the case $\tilde p_t(x_t|Y_{t-1})$ is Gaussian and $g_t$ is an affine  function of $y_t$ (and these conditions are both satisfied in the proof of Theorem \ref{th1} to characterize the saddle point), then it is not difficult to see that  there exists $\bar \lambda_t>0$ such that   $ \tilde{\Ls}$ is well defined for $\lambda_t>\bar \lambda_t$ and   $\frac{\mathrm{d}}{\mathrm{d} \lambda_t} \tilde{\Ls} (\lambda_t) \rightarrow -\infty $ as $\lambda_t \rightarrow \bar \lambda_t^+$.   Accordingly, the infimum cannot be attained for $\lambda_t \rightarrow \bar \lambda^+_t$. Since $\tilde{\Ls} (\lambda_t)$ is a continuous function  for $\lambda_t>\bar \lambda_t$, by Weierstrass' theorem it  follows that $\tilde{\Ls}$ admits a point of  minimum in $(\bar\lambda_t,\infty)$. Moreover, it is not difficult to see that $\tilde{\Ls}$ is strictly convex and thus the point of minimum is the unique stationary point. Imposing the stationarity condition we immediately see that $\lambda_t$ must satisfy condition $\Ds(\tilde{\psi}^0_{t}, {\psi}_{t})=c_t$.\qed

\end{pf}
It is worth noting that $ \tilde \psi^\star_t$ and $g^\star_t$ are mutually dependent. In particular, the minimax problem (\ref{minimax_psi}) could have more than one saddle point solution.
In what follows,  we derive the  minimizer of (\ref{minimax_psi}), for  $t\in\{0 \ldots N\},$ corresponding to the initial state $x_0$ which is Gaussian distributed.  More precisely, the idea is to substitute the structure of the maximizer $\tilde{\psi}_t^0$, as defined in Lemma 2, into Problem (\ref{minimax_psi}), obtaining an alternative formulation in which the minimizer remains unchanged. This alternative formulation allows us to characterize the minimizer.

\begin{thm}\label{th1}
Consider the estimation problem corresponding to the state space model (\ref{nomi_linea_mod_proc})-(\ref{nomi_linea_mod_obv}) and whose update estimate is obtained through (\ref{minimax_psi}).   Assume that:
\begin{equation}\label{pxY0}
p_0(x_0) \sim \mathcal{N}\left( \hat x_{0}, {P}_{0}\right).
\end{equation}
Then,  the  estimator of $x_t$ given $Y_{t}$  is
\begin{equation}\label{hatxkf}
\begin{aligned}
 \hat x_{t|t}  = \hat x_t +  L_t  (y_t -  C \hat x_t),
\end{aligned}
\end{equation}
where $L_t$ is the filtering gain: \begin{equation}\label{Lt}L_t=P_t C^\top ( C P_t C^\top +R)^{-1}.\end{equation}
Moreover, the  posterior error  covariance matrix corresponding to the least favorable posterior density is given by
\begin{equation} \label{VTT} V_{t|t} =   (P_{t|t}^{-1} -  \lambda^{-1}_t I) ^{-1}\end{equation}
where  \begin{equation}\label{ptt}
       P_{t|t}  =(I - L_t  C)P_{t}
       \end{equation}
denotes the  posterior error  covariance matrix corresponding to the  nominal posterior density   and $0<\lambda_t < \sigma_{max}(P_{t|t})$ is the unique solution to
\begin{equation} \label{thetaF}
\begin{aligned}
   \frac{1}{2}\left(\ln\det(I - \lambda^{-1}_t P_{t|t})   +  {\tr}\left((I-\lambda^{-1}_t P_{t|t})^{-1} - I\right)\right) =c_t .
\end{aligned}
\end{equation}
Moreover, \begin{equation}\label{pxY2}
\tilde p_t(x_{t+1}|Y_{t}) \sim \mathcal{N}\left( \hat x_{t+1}, {P}_{t+1}\right),
\end{equation}  the   predictor of $x_{t+1}$ given $Y_{t}$   takes the   form
\begin{equation}\label{predxkf} \hat x_{t+1}=  A \hat{x}_{t|t},\end{equation} and its corresponding prediction error covariance matrix is
\begin{equation}\label{predP}  P_{t+1} =  A V_{t|t}A^{\top}  + Q.\end{equation}
\end{thm}

\begin{pf}
We  prove (\ref{pxY2}) using the  induction principle. Condition (\ref{pxY2}) with $t=-1$ holds by (\ref{pxY0}) since there are no observations for conditioning. Next, we prove that if
\( \tilde{p}_t(x_t | Y_{t-1})  \sim \mathcal{N}\left( \hat x_{t}, {P}_{t}\right) \) for \( t \geq 0 \), then \( \tilde{p}_{t+1}(x_{t+1} | Y_t) \)  will be Gaussian.
Let \( w_t := \begin{bmatrix} x_t^{\top} & y_t^{\top} \end{bmatrix}^{\top} \). In accordance with the model for the observations (\ref{nomi_linea_mod_obv}), the nominal density of $w_t$ given $Y_{t-1}$ is \begin{equation}\nn
\begin{aligned}
  { {p}_{t}(w_{t} | Y_{t-1})} = {\psi}_{t} (y_t|x_t)\tilde  p_t\left({x}_t | Y_{t-1}\right) \sim \mathcal{N}\left(m_{t}, K_{{t}}\right)
\end{aligned}
\end{equation}  with
\begin{align*}
m_{t}&=\left[\begin{array}{c}
m_{x_{t}} \\
m_{y_t}
\end{array}\right] = \left[\begin{array}{c}
\hat x_{t}  \\
C\hat x_{t}
\end{array}\right], \\
{K}_{{t}}&=\left[\begin{array}{cc}
K_{x_{t}} & K_{x_{t}y_t} \\
K_{y_t x_{t}}  & K_{y_{t}}
\end{array}\right] =\left[\begin{array}{cc}
{P}_{{t}} &  P_t  C^\top \\
C  P_t &  C P_t C^\top + R
\end{array}\right].
\end{align*}
Then, in view of Lemma \ref{lemma1}, the least favorable density of $w_t$ given $Y_{t-1}$ takes the following form
\begin{equation}\begin{aligned} \label{tpw}
  \tilde {p}_{t}&(w_{t} | Y_{t-1}) = \tilde {\psi}^0_{t} (y_t|x_t) \tilde p_t\left({x}_t | Y_{t-1}\right)\\
  & = \frac{1}{M_{t}(\lambda_t)} \exp \left( \frac{1}{2 \lambda_t}  \left\|  x_{t}-g_t(y_t)\ \right\|^{2}\right)  p_{t}\left(w_{t} | Y_{t-1}\right).
  \end{aligned}
\end{equation}
Since $p_t(w_t|Y_{t-1})$ is Gaussian, in view of (\ref{tpw}), it follows that $\tilde p_t(w_t|Y_{t-1})$ is Gaussian.
Hence, we have
\begin{equation} \nn \begin{aligned}
\mathcal D&(\tilde{p}_{t}(w_{t} | Y_{t-1}), ~{p}_{t}(w_{t} | Y_{t-1})) \\
& :=\iint  \tilde{p}_{t}(w_{t} | Y_{t-1}) \ln\left(\frac{\tilde{p}_{t}(w_{t} | Y_{t-1})}{{p}_{t}(w_{t} | Y_{t-1})}\right)d w_t \\
& =\iint  \tilde \psi^0_t(y_t | x_t) \tilde p_t(x_t|Y_{t-1}) \ln\left(\frac{\tilde \psi^0_t}{\psi_t}\right)d y_t d x_t \\
& = \mathcal D(\tilde \psi^0, \psi_t). \end{aligned}\end{equation}
Such identity determines the ambiguity set for $\tilde{p}_{t}(w_{t} | Y_{t-1})$ corresponding to   $\mathcal{B}_t$:
\begin{equation} \nn \scalebox{0.9} {$\tilde{\mathcal{B}}_t=\{ \tilde p_t (w_{t} |  Y_{t-1}), ~  s.t.\;\;   \Ds(\tilde p_t (w_{t} |  Y_{t-1}), p_t (w_{t} |  Y_{t-1})) \leq  c_t\}.$}\end{equation}
Thus, by plugging $\tilde \psi_t^0$ into problem (\ref{minimax_psi}), we obtain the alternative minimax formulation:
\begin{equation}
 \label{minimax2psi} \underset{g_t \in \mathcal{G}_{t}}{\mathrm{min}}\max_{\tilde  {p}_{t}(w_{t} | Y_{t-1}) \in {\tilde{\mathcal{B}}}_{t}} \tilde J_t(\tilde p_{t}(w_{t} | Y_{t-1}),g_t)
\end{equation}
where the corresponding objective function is given by:
\begin{equation}\label{objt}\begin{aligned}\tilde J_t
 = \frac{1}{2} \int\left\|  x_{t}-g_{t}\left(y_{t}\right) \right\|^{2} \tilde p_{t}(w_{t} | Y_{t-1}) d w_{t}.\end{aligned}\end{equation}
In \eqref{minimax2psi} the minimizer remains the same as in \eqref{minimax_psi}, while the equivalent maximizer is $\tilde{p}_t(w_t \mid Y_{t-1})\in\tilde{\mathcal{B}}_t$.    
Since both \( {p}_{t}(w_{t} | Y_{t-1}) \) and \( \tilde{p}_{t}(w_{t} | Y_{t-1}) \) are Gaussian, then by \cite[Theorem 1]{robustleastsquaresestimation} we obtain that
the maximizer of (\ref{minimax2psi}), hereafter called $\tilde p^0_t(w_{t} | Y_{t-1})$, takes the form
\begin{equation}\label{p0wy}\tilde p^0_t(w_{t} | Y_{t-1})\sim
  \mathcal{N}(\tilde m_{t}, \tilde K_{{t}})
\end{equation} with   \begin{equation} \nn
\begin{aligned}
\tilde m_{t}=m_t =\left[\begin{array}{l}
 m_{x_{t}} \\
 m_{y_t}
\end{array}\right], \quad
\tilde {K}_{{t}}=\left[\begin{array}{cc}
\tilde K_{x_{t}} &  K_{x_{t}y_t} \\
K_{y_t x_{t}}  &  K_{y_{t}}
\end{array}\right].
\end{aligned}
\end{equation}
Moreover, 
\begin{equation}\nn
\begin{aligned}
{ \tilde{K}_{x_t}} =  V_{t|t} + L_t K_{y_t x_{t}}  L_t^\top
\end{aligned}
\end{equation}
where $L_t = K_{x_{t}y_t}  K^{-1}_{y_{t}}$,   $  V_{t|t}=( P_{t|t}^{-1}-\lambda^{-1}_t I)^{-1}$.  Notice that, the original covariance matrix $P_{t|t}$ (corresponding to the nominal posterior density) is modified by incorporating $\lambda_t^{-1}$  into $V_{t|t}$. The variable $\lambda_t$ denotes the Lagrange multiplier corresponding to the constraint $\tilde p_t(w_t|Y_{t-1})\in\tilde{\mathcal B_t}$ (and corresponding to the constraint $\tilde \psi_t\in\mathcal B_t$ in the original problem),
whose optimal value is determined by solving (\ref{thetaF}); the latter also characterizes the least favorable density in (\ref{p0wy}).  Since the model for the state process is not affected by uncertainties, we have that:
$$\tilde p_t(x_{t+1}, w_{t} | Y_{t-1}) = \tilde p^0_t(w_{t} | Y_{t-1}) p_t(x_{t+1}|x_t)$$
which is Gaussian. Thus, the corresponding marginal densities $\tilde p(x_{t+1}, y_t | Y_{t-1})$ and $\tilde p(y_t | Y_{t-1})$ are both Gaussian. Then, we obtain that $$\tilde  p_t(x_{t+1} | Y_{t}) = \tilde p_t(x_{t+1}, y_t | Y_{t-1}) / \tilde p_t(y_t | Y_{t-1}) $$
is also Gaussian.  Accordingly, we conclude that (\ref{pxY2}) holds.
Since the maximizer (\ref{p0wy}) is Gaussian, then  the corresponding minimizer  of (\ref{minimax2psi}), takes the form  in  (\ref{hatxkf})-(\ref{Lt}).
Then, (\ref{predxkf}) and (\ref{predP}) follow from the fact that $\tilde p_t(x_{t+1} | Y_{t}) $ is Gaussian.\qed
\end{pf}

\begin{algorithm}[t]
    \caption{ U-RKF at time step $t$}\label{RKF2}
    \hspace*{\algorithmicindent} \textbf{Input} $\hat{x}_{t}$, $P_{t}$, $c_t$, $y_t$,  $ A$, $ Q$, $C$, $ R$ \\
    \hspace*{\algorithmicindent} \textbf{Output}  $\hat{x}_{t|t}$, $\hat{x}_{t+1}$
    \begin{algorithmic}[1]
    \State $L_t = P_{t}  C^{\top}(CP_{t}C^{\top} + R)^{-1}$
    \State $\hat{x}_{t|t} = \hat{x}_{t}  + L_t(y_t -  C \hat{x}_{t})$
    \State $P_{t|t} = P_{t}- P_{t}  C^{\top}(CP_{t}C^{\top} + R)^{-1}C P_{t}$ \label{alg2_3}
    \State Find $ \theta_t$ s.t. $\gamma(P_{t|t}, \theta_t) = c_t$ \label{alg2_4}
    \State $V_{t|t} = (P_{t|t}^{-1} - \theta_t I) ^{-1}$ \label{alg2_5}
    \State $\hat{x}_{t+1} =  A \hat{x}_{t|t}  $
    \State $P_{t+1} = AV_{t|t}A^{\top}  + Q$
    \end{algorithmic}
\end{algorithm}

Theorem~3 essentially adopts the least favorable density within the ambiguity set characterized in (\ref{ballpsi}), which ultimately leads to  the least favorable posterior error covariance in (\ref{VTT}).
The  resulting estimator is outlined in Algorithm \ref{RKF2} where $\theta_t := \lambda^{-1}_t$ is called risk sensitivity parameter  and 
\begin{equation} \label{theta2}
\begin{aligned}
 \gamma(P, \theta) := \frac{1}{2}\left(\ln\det(I - \theta P)   +  {\tr}\left((I-\theta P)^{-1} - I\right)\right).
\end{aligned}
\end{equation}
Observe that by setting $\lambda = \lambda_t^{-1}$ and $P = P_{t|t}$ in \eqref{theta2}, we obtain the left-hand side of \eqref{thetaF}. The proposed approach differs from the classical Kalman filter in the update stage, specifically in Steps~4 and~5 of Algorithm~\ref{RKF2}, which correspond to a portion of that stage.  More precisely, the risk-sensitivity parameter $\theta_t$ is incorporated into the filtering error covariance $V_{t|t}$.
This modification effectively changes the filtering gain and makes the update stage more conservative when facing model uncertainties.  As $\theta_t$ increases, the eigenvalues of the filtering covariance matrix $V_{t|t}$ becomes more inflated, enabling the filter to hedge against potential modeling errors.  In plain words, the risk-sensitivity parameter determines how cautious the filter becomes in the presence of model uncertainty. Since this robustness mechanism is embedded solely in the update stage, we refer to our filtering approach as  update-resilient Kalman filter (U-RKF).  As we already pointed out in Section \ref{sec_2}, there exists a fundamental  difference between the proposed ambiguity set and the one in \cite{ROBUST_STATE_SPACE_LEVY_NIKOUKHAH_2013,10654520}. Such difference  is also reflected in the resulting estimators:  in U-RKF the robustification is applied to $P_{t|t}$, while in the estimators proposed in  \cite{ROBUST_STATE_SPACE_LEVY_NIKOUKHAH_2013,10654520} the robustification  is applied to $P_t$. From this viewpoint, the proposed estimator is more similar to the one in \cite{abadeh2018wasserstein} where the robustification  is applied to $P_{t|t}$ and the ambiguity set is defined with respect to  $\tilde p_t(w_t|Y_{t-1})$, i.e. the same conditional density that we have considered in the equivalent game, according to the Wasserstein distance.
However, this estimator lacks an explicit interpretation of the underlying uncertainty in terms of (\ref{nomi_linea_mod_proc}) and (\ref{nomi_linea_mod_obv}), i.e. in terms of $p_t(x_{t+1}|x_t)$ and $\psi_t(y_t|x_t)$. In particular, it has not been investigated whether this uncertainty affects only the observation model (\ref{nomi_linea_mod_obv}) or not.
 Furthermore, it lacks an  explicit characterization of the least favorable  state space model, and does not provide any filter stability  guarantee when the tolerance (i.e. Wasserstein radius) is chosen sufficiently small.
Note that, in the limiting case where $c_t = 0$ (i.e. there is no model uncertainty), both U-RKF and the one in   \cite{ROBUST_STATE_SPACE_LEVY_NIKOUKHAH_2013} coincide with the standard KF which represents the best estimator in the absence of uncertainty. 

 The least favorable model over the time interval $\{0, \ldots, N\}$ is defined in \eqref{tildepZ}, where $\tilde{\psi}_t = \tilde{\psi}_t^\star$ and $\tilde{\psi}_t^\star$ is the maximizer of \eqref{minimax_psi}. The next theorem provides a state-space representation approximating the least favorable model when the covariance matrix is taken sufficiently small in norm. The latter offers a tool for generating worst-case data within the ambiguity set. Such data can be compared with actual measurements, thereby assessing the validity of the constructed ambiguity set and enabling a systematic evaluation of different filtering algorithms under the most challenging yet admissible conditions.

\begin{thm}\label{th2} Let matrix $Q$ be taken sufficiently small in norm, the the least favorable model over the time horizon  $t\in \{0\ldots N\}$  is approximated  by:
\begin{equation}\label{LF_x2}\begin{aligned}
\eta_{t+1} &=\bar {A}_t \eta_{t} +  \bar B_t  v_t \\
y_{t} &=\bar{C}_t \eta_{t}+ \bar D_t  v_t.
\end{aligned}\end{equation}
where $\eta_t:=[\,  x_{t}^{\top}\; e_{t-1}^{\top}\;  \epsilon_{t-1}^{\top} \,]^{\top} \in \Rs^{3n}$; $ v_t :=[\,  \epsilon_{t}^{\top}\;\; \upsilon_{t}^{\top}\,]^{\top} \in \Rs^{n + m}$ is white Gaussian noise with  zero mean and covariance matrix $\Xi = \mathrm{diag}(Q, I_m)$,where $\mathrm{diag}(\cdot)$ denotes the block-diagonal matrix formed by placing its arguments along the main diagonal,  
and $\upsilon_{t} \in \Rs^m$ is   normalized white Gaussian noise. Moreover,  \begin{equation} \label{barABCD}\begin{aligned}
&\bar{A}_{t}:=\left[\begin{array}{ccc}
  A & 0 & 0\\
0 &  A- L_{t}  C A - L_t F_t A & I_n - L_{t} F_{t} - L_t C\\
0 &0 &0\\
\end{array}\right]\\
&\bar {B}_{t}:=\left[\begin{array}{cc}
 I_n  & 0 \\
 0 & -L_t  \Upsilon_t\\
 I_n & 0\\
\end{array}\right]\\
&\bar {C}_{t}:=\left[\begin{array}{llll}
C & F_t A & F_{t}
\end{array}\right], \quad \bar {D}_{t}:=\left[\begin{array}{cc}
 0 &  \Upsilon_t
\end{array}\right],
\end{aligned}\end{equation}
where   $\Upsilon_t $ is an arbitrary square root matrix  of ${O}_{t}$, i.e. $\Upsilon_t \Upsilon^\top_t = O_t$, and \begin{equation}\begin{aligned}\nn O_t&:=\left(R^{-1} -L_t^\top W_{t+1}L_t  \right)^{-1}\\
 F_{t}&:= - O_t L^\top_t W_{t+1} (I_{ n} - L_{t}  C) .\end{aligned}\end{equation}
Finally,  \begin{equation} \label{W} W_{t+1}:=   \theta_t  I_n + \Omega^{-1}_{t+1}\end{equation}
where $\Omega^{-1}_{t}$
is calculated by the following backward recursion:
\begin{equation} \label{omgea}\begin{aligned}\Omega^{-1}_{t} =  A^\top  F^\top_{t} O_t^{-1} F_{t}  A + ( A- L_{t}  C A)^\top W_{t+1} ( A- L_{t}  C A)  \end{aligned}\end{equation}
 with $\Omega^{-1}_{N+1} = 0.$
\end{thm}

\begin{pf}
Let \begin{equation}\nn e_{t}=x_{t}-\hat x_{t|t}. \end{equation}
Taking into account (\ref{nomi_linea_mod_proc}) and (\ref{nomi_linea_mod_obv}), we obtain
\begin{equation} \label{e_t1}
\begin{aligned}e_{t} = \left(  A - L_{t}  C A \right) e_{t-1} + \left( I_n - L_{t}   C \right) \epsilon_{t-1} -L_{t} \varepsilon_{t}.\end{aligned} \end{equation}
In view of  (\ref{psi_0}), with $\theta_t = \lambda^{-1}_t$, we have the least favorable density is
\begin{equation} \nn
 \tilde{\psi}_{t}^{\star} (y_t|x_t)=\frac{1}{M_{t} } \exp \left( \frac{ \theta_{t}}{2}  \left\|  e_t \right\|^{2}\right) \psi_{t}(y_t|x_t) .
\end{equation}
It is worth noting that the latter is not a normalized  density, which implies that the hostile player has the opportunity to backtrack and change the least favorable density. Since only the model (\ref{nomi_linea_mod_obv}) is affected by uncertainty, the model for the state process in (\ref{nomi_linea_mod_proc}) does not change.
Consider the nominal model  for the observations (\ref{nomi_linea_mod_obv}),
then it is not difficult to see that given $x_t$,
there is a one-to-one correspondence between $y_{t}$ and $\varepsilon_t$.   Thus, we can characterize the least favorable  model for the observations through  $\varepsilon_{t}$.
Notice that,  $\varepsilon_{t}$ does not depend on $e_{t-1}$ and $\epsilon_{t-1}$ under the nominal model. Thus, the nominal density of $\varepsilon_{t}$ is \begin{equation}\nn\varphi_t\left(\varepsilon_{t}\right)  \propto \exp \left(-\left\|\varepsilon_{t}\right\|_{R^{-1}}^{2} / 2\right)\end{equation}
where  $\propto$  means that the two terms are the same up to constant scale factors.  Instead, in view of (\ref{e_t1}),  we make the guess that  the least favorable density of $\varepsilon_{t}$  is related  to $e_{t-1}$ and $\epsilon_{t-1}$, namely, we consider $\tilde \varphi_t(\varepsilon_{t} | e_{t-1}, \epsilon_{t-1}).$   Accordingly, we construct the term
\begin{equation} \nn \exp\left( \frac{1}{2} \left\| e_{t}\right\|_{ \Omega_{t+1}^{-1}}^{2}   \right)\end{equation}
to indicate the cumulative error  of the retroactive probability density changes of $\varepsilon$ over the interval $[0,t]$. Here, $\Omega_{t}$ is a positive definite matrix of dimension $n$.
 Thus, the least favorable density of $\varepsilon$ over the time interval  $ \{t+1\ldots N\}$ takes the form:
\begin{equation}\label{lfd_t}\begin{aligned}
 \prod_{s=t+1}^{ N} &\exp \left(\frac{ \theta_{s} }{2} \left\| e_{s} \right\|^{2} \right) \varphi_s\left(\varepsilon_{s}\right)\\
& \propto \exp\left( \frac{1}{2}\left\| e_{t}\right\|_{ \Omega_{t+1}^{-1}}^{2} \right)\prod_{s=t+1}^{ N} \tilde {\varphi}_{s}\left(\varepsilon_{s} |  e_{s-1}, \epsilon_{s-1}\right).
\end{aligned}\end{equation}
Notice that,  if the matrix $\Omega_{t}$ can be evaluated recursively, it is then possible to find the least favorable density $\tilde \varphi_t(\varepsilon_{t} | e_{t-1}, \epsilon_{t-1})$  through a backward recursion. Decreasing the time index $t$ by 1 in (\ref{lfd_t}) and
subtracting it in (\ref{lfd_t}), it is not difficult  to see that
\begin{equation*}\begin{aligned}
\tilde \varphi_{t}& \left(\varepsilon_{t}| e_{t-1}, \epsilon_{t-1}\right)  \propto \nn\\ 
&\exp  \left(-\frac{1}{2} \left( -\left\| e_{t}\right\|^{2}_{W_{t+1}} + \left\| e_{t-1}\right\|^{2}_{ \Omega^{-1}_{t}}
 + \left\|\varepsilon_t\right\|^{2}_{R^{-1}} \right)\right),
\end{aligned}\end{equation*}
where $W_{t+1}$  is defined in (\ref{W}). Moreover, taking into account (\ref{e_t1}), we obtain the expression in (\ref{varphi_2}) {(see on the top of the next page)}.
\begin{figure*}
\begin{equation} \label{varphi_2}
\begin{aligned}
\tilde \varphi_{t} \left(\varepsilon_{t}| e_{t-1}, \epsilon_{t-1}\right)
 \propto \exp & \left(-\frac{1}{2} \left( e^\top_{t-1} ( \Omega^{-1}_{t} -( A- L_{t}  C A)^\top W_{t+1} (A- L_{t}  C A) ) e_{t-1}
\right.\right. \\
 &   - \epsilon^\top_{t-1}   (I_n - L_{t}   C)^\top W_{t+1}  (I_n - L_{t}   C) \epsilon_{t-1} + \varepsilon^\top_t (R^{-1} -L_t^\top W_{t+1}L_t  ) \varepsilon_t\\
 & + 2  e^\top_{t-1} ( A- L_{t}  C A)^\top W_{t+1}  L_t \varepsilon_t +
 2 \epsilon^\top_{t-1} (I_n - L_{t}  C)^\top W_{t+1}    L_t \varepsilon_t \\
 & \left.\left. - 2  e^\top_{t-1} ( A- L_{t}  C A)^\top W_{t+1} (I_n - L_{t}   C)\epsilon_{t-1} \right)\right).
\end{aligned}\end{equation}
	{\noindent} \rule[-10pt]{18cm}{0.05em}
\end{figure*}
Assuming that the covariance matrix $Q$ is small in norm,  $\epsilon_{t-1}$ is likely to take small values; therefore, the terms $ \epsilon^\top_{t-1}   (I_n - L_{t}   C)^\top W_{t+1}  (I_n - L_{t}   C) \epsilon_{t-1}$ and $   e^\top_{t-1} ( A- L_{t}  C A)^\top W_{t+1} (I_n - L_{t}   C)\epsilon_{t-1}$ can be neglected.  Under such approximation $\tilde\varphi_t$  can be expressed in the following compact way:
\begin{equation} \label{lef_v2}  \begin{aligned}
\tilde \varphi_{t}& \left(\varepsilon_{t}| e_{t-1}, \epsilon_{t-1}\right) \\&  \propto
\exp  \left(-\frac{1}{2}   \left\| \varepsilon_t - (F_{t} A_{t-1}  e_{t-1} + F_{t}  \epsilon_{t-1})\right\|_{ O_t^{-1}}^{2}   \right)
 \end{aligned}\end{equation}
where the backward recursion (\ref{omgea}) is obtained by matching the quadratic term of $e_{t-1}$ in (\ref{varphi_2}) with the one in (\ref{lef_v2}). In view of (\ref{lef_v2}), we have
\begin{equation}\nn \tilde \varphi_{t} \left(\varepsilon_{t}| e_{t-1}, \epsilon_{t-1}\right) \sim \mathcal{N}\left(F_{t}  A e_{t-1} + F_{t}  \epsilon_{t-1}~, ~ O_t\right).\end{equation}
Thus, the guess   that  $\varepsilon_{t}$   depends on $e_{t-1}$ and $\epsilon_{t-1}$  holds. We conclude that
\begin{equation}\label{varpsilon}\varepsilon_t= F_{t}  A e_{t-1} + F_{t}  \epsilon_{t-1} + \Upsilon_t \upsilon_{t}\end{equation}
where $\upsilon_{t}$ is normalized white Gaussian noise. Substituting (\ref{varpsilon}) in (\ref{nomi_linea_mod_obv}) and (\ref{e_t1}), we obtain:
\begin{equation}\nn \begin{aligned}
e_{t} &= \left(  A - L_{t}  C A \right) e_{t-1}+ \left( I_n - L_{t}   C \right) \epsilon_{t-1}  \\
&\hspace{0.5cm} -L_{t} (F_{t}  A   e_{t-1} + F_{t}  \epsilon_{t-1} + \Upsilon_t \upsilon_{t}) \\
y_{t} &= C x_{t} + F_{t}   A   e_{t-1} + F_{t}  \epsilon_{t-1} + \Upsilon_t \upsilon_{t}.
\end{aligned}\end{equation}
Taking into account (\ref{nomi_linea_mod_proc}), $\eta_t:=[\,  x_{t}^{\top}\; e_{t-1}^{\top}\;  \epsilon_{t-1}^{\top} \,]^{\top} $, $ v_t :=[\,  \epsilon_{t}^{\top}\;\; \upsilon_{t}^{\top}\,]^{\top} $, and $\bar A_t,~\bar B_t, ~\bar C_t, \bar D_t $ as in Eq. (\ref{barABCD}), then we obtain the following least favorable state space model:
\begin{equation*}\begin{aligned}
\eta_{t+1} &=\bar {A}_t \eta_{t} +  \bar B_t  v_t \\
y_{t} &=\bar{C}_t \eta_{t}+ \bar D_t  v_t,
\end{aligned}\end{equation*}
which corresponds to the final form  in (\ref{LF_x2}).   \qed
\end{pf}

It is interesting to point out that the proposed ambiguity set \eqref{ballpsi} leads to a least favorable model whose structure is different  from the one  corresponding to \eqref{oldAB}, see \cite[Section V]{ROBUST_STATE_SPACE_LEVY_NIKOUKHAH_2013}. In particular, the former has state space dimension equal to $3n$, while the one in \cite{ROBUST_STATE_SPACE_LEVY_NIKOUKHAH_2013} has dimension $2n$. We conclude that, uncertainty only in the observation model leads to a least favorable model which is more complex than the one where uncertainty is in both the state process and observation models.

The theorem above also suggests the operative way to compute the least favorable model: first, we perform a forward recursion to compute $L_t$ through Algorithm \ref{RKF2}; then, we perform a backward recursion to compute $\Omega^{-1}_{t}$; finally, we compute the state space matrices of (\ref{LF_x2}).

We conclude  this section  by showing how to evaluate the performance of an arbitrary state   estimator of the form
\begin{equation} \label{x_prime}\hat x^\prime_{t|t} =  A \hat{x}^\prime_{t-1|t-1} + L^\prime_t(y_t -  CA \hat{x}_{t-1|t-1})\end{equation} where   $L^\prime_t$ is a Kalman gain sequence. Notice that, if we take $L^\prime_t = L_t$, we obtain the U-RKF.
Let $e^\prime_t = x_t - \hat x^{\prime}_{t|t}$ be the estimation error corresponding to (\ref{x_prime}), and let $$\begin{aligned}
  \Delta^{\prime}_t &  =  A - L_{t}^{\prime}  C A; ~  ~  \Delta_t  =  A - L_{t} C A ; \\
\Lambda^{\prime}_t &= I_n - L^{\prime}_{t} F_{t} - L^{\prime}_t  C; ~ \Lambda_t = I_n - L_{t} F_{t} - L_t  C.
\end{aligned}$$ Then, the covariance matrix of   $ \bar e_t:=[\, (e^\prime_{t})^{\top} \; e_{t}^{\top} \; \epsilon_{t}^{\top} \,]^{\top}$ is obtained through the   Lyapunov equation
\al{\label{a_lfm1}\Pi_{t+1}=\Gamma_t \Pi_t \Gamma_t^\top +  \text{X}_t \Xi \text{X}_t^\top}
where
{\small\begin{equation}  \begin{aligned}
\Gamma_{t} =  \left[\begin{array}{ccc}
\Delta^{\prime}_t &  - L_{t}^{\prime} F_{t} A & \Lambda^{\prime}_t  \\
 0 & \Delta_t - L_t F_{t}  A & \Lambda_t \\
 0 & 0 & 0
\end{array}\right],\;
{\text X}_t=\left[\begin{array}{cc}
0  & -L^{\prime}_t  \Upsilon_t\\
0  & -L_t  \Upsilon_t \\
I_n & 0
\end{array}\right] .
\end{aligned}\nn
\end{equation}}Thus, the covariance matrix of the estimation error $ e^\prime_t$ is given by the $n \times n$ submatrix of $\Pi_{t}$ in the top-left position.

\section{Filter stability} \label{sec_4}
In this section we consider the situation  in which  the tolerance is constant, i.e. $c_t=c$.   We assume that the pair $(A,C)$ is observable. Notice that, the pair $(A,Q)$ is reachable because $Q>0$.
In what follows, we will show that it is possible to characterize an upper bound for the tolerance, say $c_{MAX}$, which guarantees that the gain $L_t$ of U-RKF converges to a constant value for any $c\in (0,c_{MAX}]$. Then, we will show that for $c>0$ taken sufficiently small it is possible to guarantee in steady state that the estimation error under the least favorable model is bounded in mean square and the least favorable model is a state space model with constant parameters.

In view of Algorithm \ref{RKF2}, $P_t$ is  characterized
through the recursion
\al{\label{rec} P_{t+1}=r_c(P_t):=A(P_t^{-1}+C^\top R^{-1}C-\theta_t I)^{-1}A^\top +Q}
where we exploited the fact that the equation in Step \ref{alg2_3} can be written as
$P_{t|t}=(P_t^{-1}+C^\top R^{-1}C)^{-1}.$
It is worth noting that the ``distorted'' Riccati operator $r_c$ has the same structure of the one  for   the prediction-resilient Kalman filter (P-RKF) proposed in \cite{ROBUST_STATE_SPACE_LEVY_NIKOUKHAH_2013}: the difference regards how $\theta_t$ depends on $P_t$. Such difference requires a convergence analysis which is substantially different to the one for the prediction-resilient case.
Let $k\geq n$, we define the matrix \al{
\bf R_k & :=  {\cal O}_k^{\top}({\cal Q}_k +{\cal H}_k{\cal H}_k^{\top})^{-1}{\cal O}_k + {\cal J}_k^{\top} S_k^{-1} {\cal J}_k \nn}
where
\al{S_k &:= {\cal L}_k(I+{\cal H}_k^{\top} {\cal Q}_k^{-1}{\cal
H}_k)^{-1}{\cal L}_k^{\top} - \phi_k^{-1} \otimes I \nn\\
 {\cal J}_k& :=  {\cal O}_k^{R}-{\cal L}_k {\cal H}_k^{\top}[{\cal Q}_k +{\cal H}_k {\cal H}_k^{\top}]^{-1} {\cal O}_k \nn  \\  {\cal O}_k & :=[\;
 (CA^{N-1})^{\top} \; \ldots \;  (C A)^{\top} \;  C^{\top}\;]^{\top} \nn \\
{\cal O}_k^R & :=  [\;
 (A^{N-1})^{\top} \; \ldots \;  A^{\top} \; I
\;]^{\top}  \nn\\
{\cal Q}_k & := I \otimes  R\nn\\
{\cal H}_k & :=\mathrm{Tp} \left(
0 \; H_1,H_2 , \ldots, H_{N-2} , H_{N-1}
\right)\nn\\
{\cal L}_k & :=\mathrm{Tp} \left(
0  , L_1 , L_2  , \ldots  , L_{N-2}  ,L_{N-1}
\right)\nn\\
H_t & := C A^{t-1} Q^{1/2},\;L_t  :=
A^{t-1} Q^{1/2} \nn}
where $\mathrm{Tp}(\cdot)$ denotes the block upper triangular Toeplitz matrix whose argument define its first block row and $Q^{1/2}$ denotes a square root matrix of $Q$. Matrix $\mathbf R_k$ is related to the
observability Gramian of the $k$-fold composition of the mapping $r_c(\cdot)$, see \cite[Section 4]{LEVY_ZORZI_RISK_CONTRACTION} for more details  about mappings of this form.
Let $\phi_k\in(0,\sigma_{max}({\cal L}_k(I+{\cal H}_k^{\top} {\cal Q}_k^{-1}{\cal
H}_k)^{-1}{\cal L}_k^{\top}))$ be the maximum value for which $\mathbf R_k$ is a positive definite matrix. As explained in \cite{LEVY_ZORZI_RISK_CONTRACTION}, such $\phi_k$ does exist because $(A,C)$ is observable. By Proposition 3.1 in \cite{zorzi2017convergence}, if \al{\label{cond_th}\theta_t\leq \phi_k, \; \forall t\geq q+1}
for some $q\in \mathbb N$, $(A,Q)$ is reachable and $(A,C)$ is observable,  then the sequence generated by (\ref{rec}) converges and the corresponding algebraic equation admits a unique solution $P>0$.  So, we have to find the condition on $c$ for which (\ref{cond_th}) holds. In what follows, we consider the sequence generated by the Riccati operator
\al{\ \bar P_{t+1}&=\bar r(\bar P_t):=A(\bar P_t^{-1}+{ C^\top}R^{-1}C)^{-1}{ A^\top}+Q,\nn\\
 \bar P_0&=Q\nn}
 and we define $\bar P_{t|t}=(\bar P_t^{-1}+{C^\top}R^{-1} C)^{-1}$.

\begin{thm} \label{teo_c} Let model (\ref{nomi_linea_mod_proc})-(\ref{nomi_linea_mod_obv}) be such that $(A,C)$
is observable. We define
\al{c_{MAX}=\gamma(\bar P_{q|q},\phi_k )>0 \nn}
where $k\geq n$ and  $q\in\mathbb N$. Let $c$ be such that $c\in(0,c_{MAX}]$, then  the sequence generated by the iteration (\ref{rec}) converges to a unique solution $P>0$ for any $P_0>0$. Moreover, the limit $L$ of the filtering gain $L_t$ as $t\rightarrow \infty$ has the property that $A(I-LC)$
 is stable.
\end{thm}
\begin{pf} By Lemma 4.1 in \cite{zorzi2017convergence}, we have that
\al{P_t\geq \bar P_q,\quad t\geq q+1\nn}
for any $q\geq 0$. Thus,
\al{(P_t^{-1}+C^\top R^{-1}C)^{-1}\geq (\bar P_q^{-1}+{ C^\top}R^{-1}C)^{-1} \nn} which  implies
\al{P_{t|t}\geq \bar P_{q|q}, \quad t\geq q+1.\nn }
Next, we prove that (\ref{cond_th}) holds: we prove that by contradiction, i.e. we assume that $\theta_{t}>\phi_k$ for some $t\geq q+1$.
Since $\gamma(X,\cdot)$, with $X\geq0$ and $X\neq 0$, is monotone increasing over $\mathbb R_+$ and $\gamma(X,\theta)\geq \gamma(Y,\theta)$ for $X\geq Y$, see  Lemma 4.3 in \cite{zorzi2017convergence}, it follows that
\al{c=\gamma(P_{t|t},\theta_t)\geq  \gamma(\bar P_{q|q},\theta_t)>\gamma(\bar P_{q|q},\phi_k)=c_{MAX}\nn}
which is a contradiction. Accordingly, all the hypotheses of Proposition 3.1 in \cite{LEVY_ZORZI_RISK_CONTRACTION} hold and thus  (\ref{rec}) converges to a unique solution $P$. The stability of $A(I-LC)$ follows from the fact that the algebraic equation corresponding to (\ref{rec}) and having unique solution $P$, can be written as the Lyapunov equation
\al{P=A(I-LC)P(I-LC)A^\top+Q+AL(AL)^\top.\nn}\qed
\end{pf}

The least favorable model has been characterized over a simulation horizon  $ \{0\ldots N\}$ through  a forward and a backward recursion.  Then, the least favorable model in steady state  is obtained in the interval  $t\in \{ \lfloor \alpha N\rfloor \ldots \lceil \beta N\rceil\}$  as $N\rightarrow \infty$ where $\alpha$ and $\beta$ are such that $0<\alpha<\beta<1$.
Let  $t\in \{ \lfloor \alpha N\rfloor \ldots  N\}$. By Theorem \ref{teo_c},  $P_t\rightarrow P$, $\theta_t\rightarrow \theta$ and $L_t\rightarrow L$  as $ N\rightarrow \infty$.
  So, the backward recursion (\ref{omgea}) becomes
\al{\label{back_steady}\Omega_t^{-1}=(F_tA)^\top O_t^{-1} F_tA+ \bar A^\top W_{t+1}\bar A}
with
\al{&O_t=( I-L^\top W_{t+1} L)^{-1}, & &F_t=-O_tL^\top W_{t+1}  (I-LC)\nn\\
&\bar A= (I-LC)A, & & W_{t}= \Omega_{t}^{-1} +\theta I.\nn}
If $\Omega_t^{-1}$ (or equivalently  $W_{t}$) converges as $ N\rightarrow \infty$, then matrices $\bar A_t$, $\bar B_t$, $\bar C_t$, $\bar D_t$ converges and thus the least favorable model is a state space model with constant parameters in the steady state interval  $t\in \{ \lfloor \alpha N\rfloor \ldots   N \}$  as $ N\rightarrow \infty$. Notice that, $\bar A$ and $L$ depend on $c$ through $\theta$.
\begin{prop}  \label{teoLFConv}\label{teoLFConv}Assume that the map
\al{f \, :\, [0,\check \theta]& \rightarrow   \mathbb R^{n\times n}\times \mathbb R^{n\times m}\nn\\
  \theta  & \mapsto   (\bar A, L)\nn}
  is continuous for $\check \theta>0$ sufficiently small. Then, there exists $c>0$ sufficiently small such that $W_t$,  with  $t\in \{  0 \ldots \lceil \beta N\rceil\}$,  converges to $W$ as $ N\rightarrow \infty$. Moreover, let $J=\bar A^\top WL({ L^\top}W L-I)^{-1}$, then $\bar A- L{J^\top}$ is a stable matrix.
\end{prop}
\begin{pf} First, notice that (\ref{back_steady}) can be written as
\al{\Omega_t^{-1}&=\bar A^\top[W_{t+1}+W_{t+1}L(I-L^\top W_{t+1}L)^{-1}L^\top W_{t+1}]\bar A\nn\\
&=\bar A^\top (W_{t+1}^{-1}-LL^\top)^{-1}\bar A.\nn}
Adding on both sides $\theta I$ we obtain
\al{W_{t}= \bar A^\top(W_{t+1}^{-1}-LL^\top)^{-1}\bar A+\theta I \nn}
which is a Riccati recursion with terminal condition $W_{N}=\theta I$. The latter is similar to the one  considered in \cite{zorzi2018robust} and using similar reasonings it is possible to prove that  if $\theta>0$ is sufficiently small  (and thus $f$ is continuous in a neighborhood of $\theta$), then $W_t$ converges and $\bar A- LJ^{\top}$  is a stable matrix.
Notice that, $P$ and $\theta$ are related through
\al{\label{ineqPQ0} c=\gamma((P^{-1}+C^\top R^{-1}C)^{-1},\theta).}
Moreover, $P$ solves the algebraic form of the Riccati recursion (\ref{rec}), so $P\geq Q$ and thus \al{\label{ineqPQ}(P^{-1}+C^\top R^{-1}C)^{-1}\geq (Q^{-1}+C^\top R^{-1}C)^{-1}>0.}
Recall that $\gamma(X,\cdot)$, with $X\geq 0$ and $X\neq 0$, is monotone increasing over $\mathbb R_+$ and $\gamma(X,\theta)\geq \gamma(Y,\theta)$ for $X\geq Y$, \cite[Lemma 4.3]{zorzi2017convergence}. Moreover, $\gamma(X,0)=0$ for any $X\geq0$ and
the range of $[0,\sigma_{max}(X)^{-1})$  under $\gamma(X,\cdot)$ is $[0,\infty)$. Taking into account (\ref{ineqPQ0})-(\ref{ineqPQ}), we conclude that it is possible to take $c$ sufficiently small such that $\theta$ is arbitrary small.\qed
\end{pf}

Let $e_t=x_t-\hat x_{t|t}$ be the state estimation error of U-RKF. The covariance matrix of $e_t$  is given by the $n\times n$ top-left submatrix of $\Pi_t$ which  obeys the recursion in (\ref{a_lfm1}) with $L_t^\prime=L_t$. In the case the least favorable model is in steady state, then the recursion takes the form  as in (\ref{a_lfm1}). Since $\Gamma_t$ and ${\text X}_t$ depend on the constant parameters of the least favorable model and the filtering gain of U-RKF, we have $\Gamma_t\rightarrow \Gamma$, $ \text{X}_t\rightarrow \text X$ as $t\rightarrow\infty$ and {\small\al{\Gamma= \left[\begin{array}{ccc}
(I-LC)A &    - L FA & I_n - L F - L C \\
0 & (I- LC)A-LFA\;\;\; & I_n - L F - L C \\
0 & 0 & 0
\end{array}\right].\nn}
Matrix $\Gamma$ is a stable matrix because its eigenvalues are the ones of $(I-LC)A$ and $(I- LC)A-LFA$ plus the ones in the origin.} By  Theorem \ref{teo_c}, $A(I-LC)$ is a stable matrix for $c$ sufficiently small. Since $A(I-LC)$ and $(I-LC)A$ have the same nonnull eigenvalues, we conclude that $(I-LC)A$ is stable.  Moreover,
$$(I- LC)A-LFA= \bar A- JL^\top$$
which is a stable matrix for $c$ sufficiently small by Proposition~\ref{teoLFConv}. Since $\Gamma_t$ converges to a stable matrix, by \cite[Theorem 1]{cattivelli2010} we have that $\Pi_t$ converges to the unique solution of the algebraic form of the Lyapunov recursion in (\ref{a_lfm1}). Accordingly, the covariance matrix of $e_t$ is bounded in steady state.


\section{Numerical examples}\label{sec_ne}
We present some numerical results to  assess the proposed estimator. First, we analyze the performance of the filter in the worst case scenario. Then, we test it in a mass-spring-damper (MSD) system where uncertainties are primarily concentrated in the measurements.

\subsection{Worst case analysis}
\begin{figure}[t]
  \centering
  \includegraphics[width=3in]{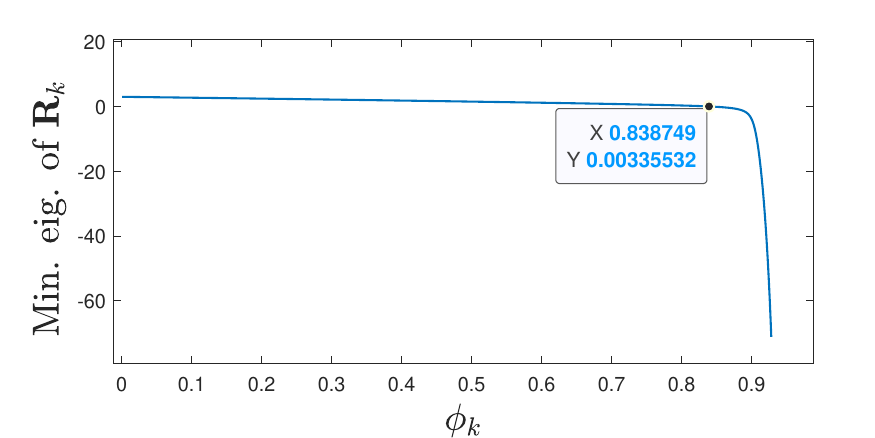}
  \caption{ Minimum eigenvalue of $\mathbf R_k$ as a function of $\phi_k$ with $k=10$. The largest value of $\phi_k$ such that $\mathbf R_k$ is positive definite is approximately equal {\mz to 0.839}.}\label{F_1}
\end{figure}

We consider a nominal state space model of the form (\ref{nomi_linea_mod_proc})-(\ref{nomi_linea_mod_obv}) with constant matrices
\begin{equation}\label{model1} \mz
\begin{aligned}
& A=\left[\begin{array}{cc}
0.95 & 1 \\
0 & 1.2
\end{array}\right], \quad Q=\left[\begin{array}{ccc}
    0.01   & 0 \\
     0 &  0.01
\end{array}\right], \\
& C=\left[\begin{array}{cc}
1 & 0 \\
0.1 & -0.1
\end{array}\right], \quad  R= \left[\begin{array}{ccc}
    1   & 0 \\
     0 &  1
\end{array}\right],
\end{aligned}
\end{equation}
and the initial state is modelled as a Gaussian random  vector with zero mean and covariance matrix $\mz P_0=I$.
Since such model is both reachable and observable, we compute the upper bound $c_{MAX}$ for which we know that the proposed U-RKF converges.
Given the matrices in (\ref{model1}), setting $k=10$ and $q=20$,  then we found $\phi_k$ is approximately equal {\mz to 0.837}, as shown in Fig. \ref{F_1}, and $$\bar P_{q|q} = \mz \left[\begin{array}{cc}
0.4609 & 0.1727 \\
0.1727 &  0.1036
\end{array}\right].$$
Thus, by Theorem \ref{teo_c}, we have that $\mz c_{MAX} = 0.5001$.

\begin{figure}
  \centering
  \includegraphics[width=3in]{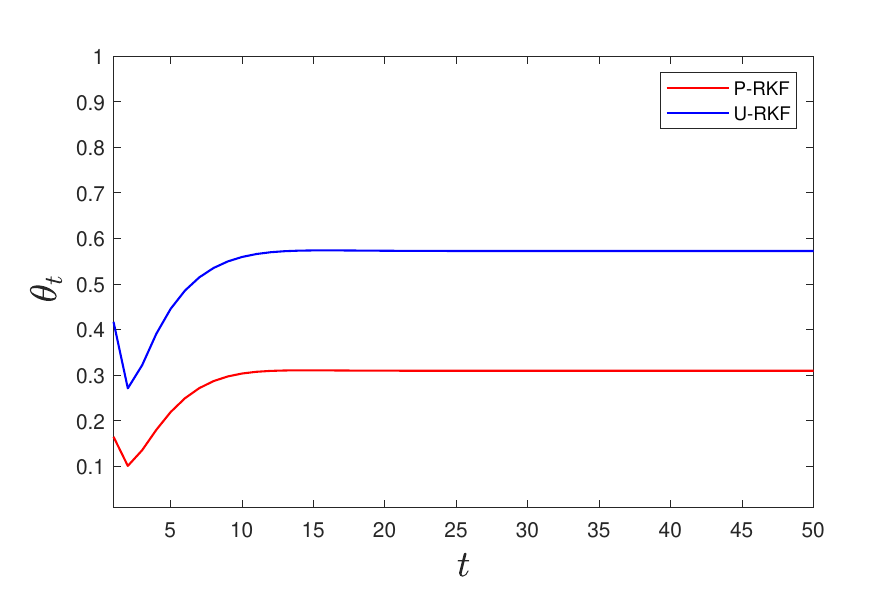}
  \caption{Risk sensitivity parameter $\theta_t$ of  P-RKF (red line)   and U-RKF (blue line) when $\mz c=0.1$.}\label{F_4}
\end{figure}

Next, we compare the performance of these estimators  when applied to the least favorable model derived in Section \ref{sec_3}.  The constant tolerance, such that $c\leq c_{MAX}$, is the same for U-RKF, P-RKF and the least favorable model.
\begin{figure}[t]
  \centering
  \includegraphics[width=3in]{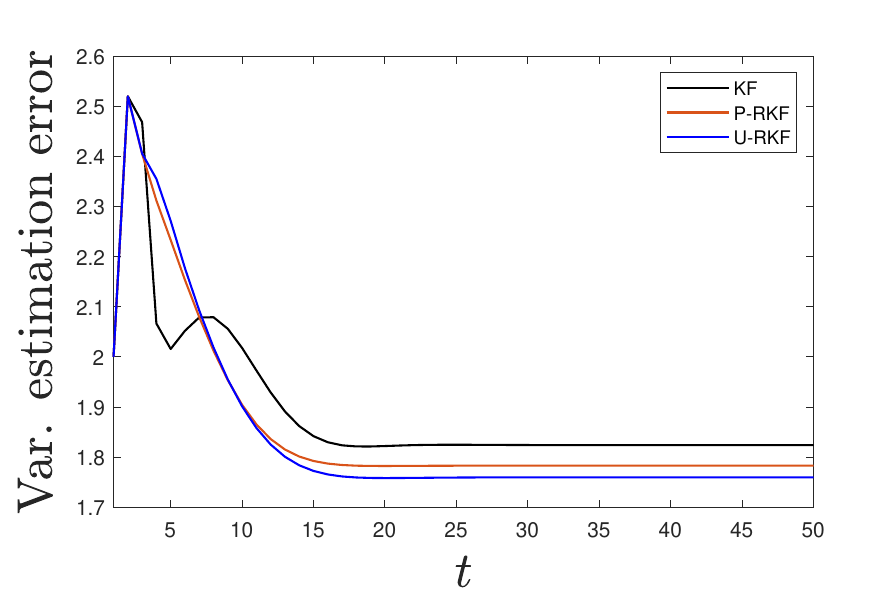}
  \caption{ Variance  of the estimation error when KF (black line), P-RKF (red line) and U-RKF (blue line)  are applied to the least favorable model with $\mz c= 0.1$.}\label{F_2}
\end{figure}
\begin{figure}[t]
  \centering
  \includegraphics[width=3in]{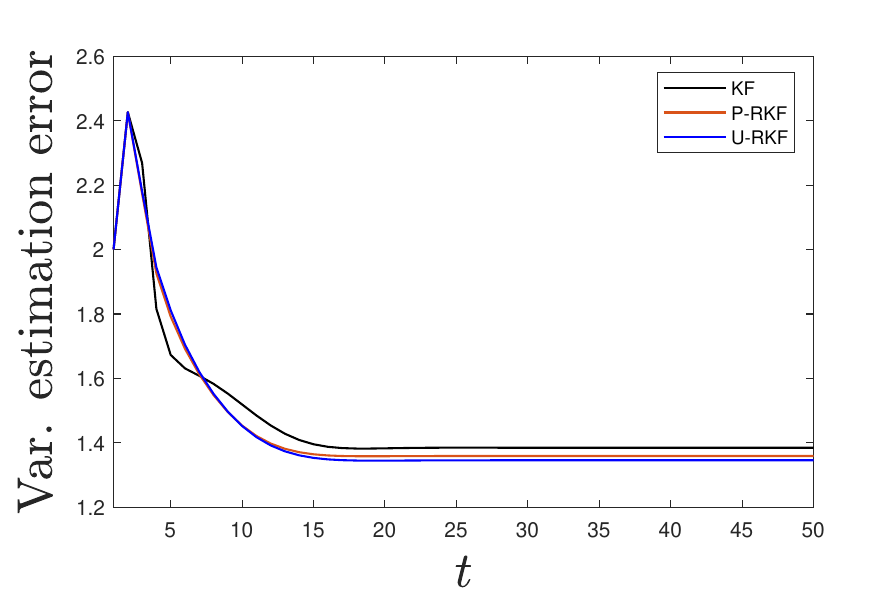}
  \caption{ Variance of the estimation error when  KF (black line), P-RKF (red line) and U-RKF (blue line)  are applied to the least favorable model with $\mz c= 5 \cdot 10^{-2}$.}\label{F_3}
\end{figure}Fig. \ref{F_2} and Fig. \ref{F_3} show  the  variance of the estimation error of the estimators,  computed through (\ref{a_lfm1}) for {\mz $c= 0.1$ and $c= 5\cdot 10^{-2}$}.
As we can see, the error variance converges for all the estimators. Also, the proposed U-RKF  outperforms both P-RKF  and KF.
Moreover, two interesting and relevant aspects emerge. First, the larger the tolerance \( c \) is, the more evident the difference in performance becomes.  Second, P-RKF outperforms KF. In plain words,  P-RKF accounts for model uncertainties, but it is overly risk adverse in the sense that it considers model uncertainties also in the process equation (\ref{nomi_linea_mod_proc}).
Finally, Fig. \ref{F_4} illustrates the risk sensitivity parameter for U-RKF  and P-RKF  with  \(\mz  c = 0.1 \).  As we can see, although their Riccati recursions have the same structure, their risk sensitivity parameters  are different.
Finally, we also considered the situation in which the least favorable model is obtained by introducing small perturbations in the matrices $A$ and $Q$ (i.e. in the state process model). We found that U-RKF still outperforms both P-RKF and the KF.

\subsection{MSD under sensor uncertainties}\label{sec_ne2}
\begin{figure}[t]
\centering
\begin{tikzpicture}[scale=0.8, thick]

    \draw[thick] (-2,-1.5) -- (2,-1.5);

    \draw[fill=blue!20] (-1.5,0) rectangle (1.5,1);
    \node at (0,0.5) {\Large ${m}$};

    \draw[thick] (-1.2,0) -- (-1.2,-0.2);
    \draw[thick, gray] (-1.2,-0.2) -- (-1.4,-0.3) -- (-1.0,-0.4)
                        -- (-1.4,-0.5) -- (-1.0,-0.6) -- (-1.4,-0.7)
                        -- (-1.0,-0.8) -- (-1.4,-0.9) -- (-1.0,-1.0)
                        -- (-1.2,-1.1);
    \draw[thick] (-1.2,-1.1) -- (-1.2,-1.5);
    \node[left] at (-1.5,-0.7) {\Large ${k}$};  

    \draw[thick, gray] (1.1,0) -- (1.1,-0.6);  
    \draw[thick, gray] (0.7,-0.6) rectangle (1.5,-1.2); 
    \draw[thick, gray] (1.1,-1.2) -- (1.1,-1.4);  
    \draw[thick] (1.1,-1.4) -- (1.1,-1.5);  
    \node at (1.1,-0.9) {\Large $\bf{c}$};  

    \draw[thick, ->] (2,0) -- (2,0.5) node[right] {\Large $p$};
    \draw[thick] (1.8,-0) -- (2.2,-0);

\end{tikzpicture}
\caption{Mass-spring-damper system.}  
\label{fig:mass_spring_damper}  
\end{figure}
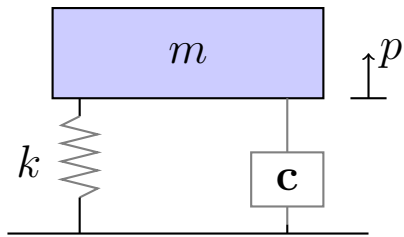
We consider a mass-spring-damper system, as shown in Fig. \ref{fig:mass_spring_damper}. The equation of the motion for this system is given by:
\begin{equation}\label{msdm}
m \ddot{p}   +  {\bf c} (\dot{p}+\nu) +k p= F
\end{equation}
where $p$ represents the displacement of the object with mass $m = 0.1~ (kg)$  away from its resting position, $k = 5~ (N/m)$ is the spring constant, ${\bf c} = 2~ (Ns/m)$ is the damping coefficient and $F$ denotes the external force, which is white Gaussian  noise with zero mean and variance equal to 0.9; $\nu$ is white Gaussian  noise with zero mean and variance equal to 0.09 which  corresponds to the presence of a small ``disturbance'' force acting on the damper (e.g. the force generated by road irregularities, such as bumps and surface roughness, in a car suspension system).
The displacement $p$ is measured using a sensor with sampling time $T_s = 0.1s$. 
 Let $y_t$ denote the sensor measurement at time $t$,  we consider  the following four classical types of sensor uncertainties:
\begin{itemize}
  \item Sensor drift $$y_t = p_t + \tilde \varepsilon_t, ~ \tilde \varepsilon_t \sim \mathcal N(0.1, R);$$
  \item Uniform noise $$y_t = p_t + \tilde \varepsilon_t, ~ \tilde \varepsilon_t \sim \mathcal U(-0.9, 1.1);$$
  \item Nonlinearity (dead zone) $$y_t = \left\{\begin{array}{ll} p_t +  \tilde \varepsilon_t & \text{when} ~ |p_t + \tilde \varepsilon_t| \geq 0.1 \\ 0 & \text{when}~ |p_t + \tilde \varepsilon_t| < 0.1  \end{array}\right. $$
where  $\tilde \varepsilon_t \sim \mathcal N(0, R)$;
        \item Outlier-contaminated noise  $$y_t = p_t + \tilde \varepsilon_t, ~\tilde \varepsilon_t \sim\left\{\begin{array}{lll} \mathcal N\left({0}, R\right) & \text { w.p. } & 0.9 \\ \mathcal N\left({0}, 5 R\right) & \text { w.p. } & 0.1\end{array}\right. $$ where
 w.p. stands for ``with probability''.
\end{itemize}
In the cases above $R = 0.25$. For each type of sensor uncertainty, we generate $M=1000$ state and measurement trajectories with $N=200$, i.e. the total time is $20s$.

Our aim is to estimate the displacement using the sensors measurements previously generated. In doing that, we consider the state space model of (\ref{msdm}) with state $x = [\,p\; \;\dot{p}\,]^\top$ where we neglect the presence of the small disturbance force $\nu$. Then, we discretize it with sampling time  $T_s$ obtaining an equation of the form (\ref{nomi_linea_mod_proc}). Finally, we assume that the nominal force is white Gaussian noise with zero mean and unit variance, i.e. its variance is slightly different than the actual one. In plain words, there is a mild mismatch between the nominal and  the actual state process model.  The nominal observation model is equal to (\ref{nomi_linea_mod_obv}) with $C=[\,1 \;\; 0\, ]^\top$ and $\varepsilon _t$ is white Gaussian noise with zero and variance $R$. Then, we set  $\hat x_0 = [\, 0\;\;0\,]^\top$ and  $P_0 = 0.05I$.

\begin{figure}[t]
  \centering
  \includegraphics[width=3.2in]{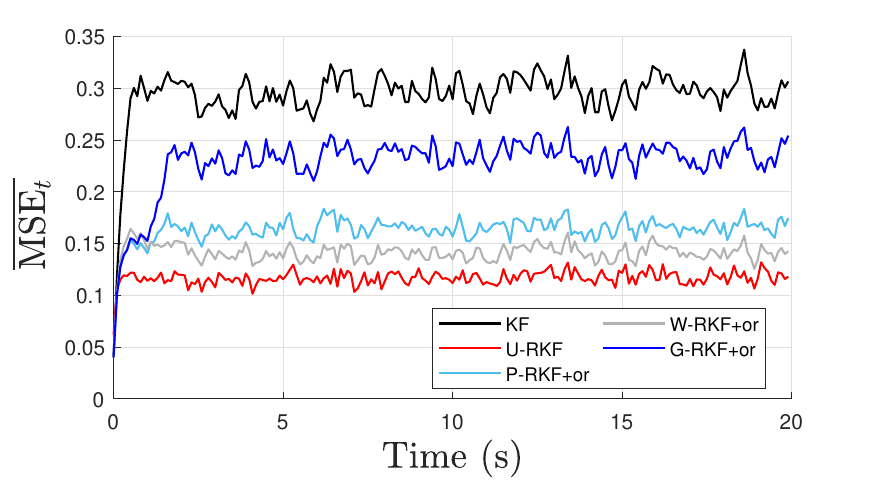}
  \caption{Average variance of the displacement error for the different filters in the presence of  sensor drift.}\label{F_drif}
\end{figure}
\begin{figure}[t]
  \centering
  \includegraphics[width=3.2in]{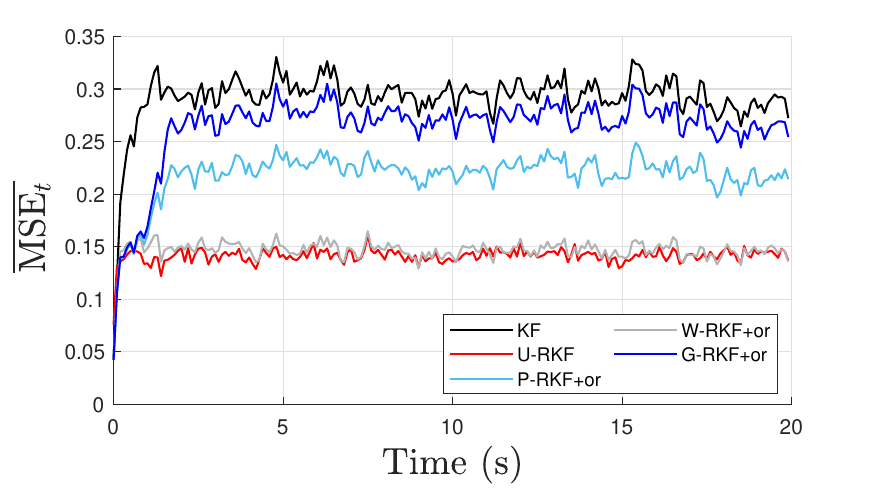}
  \caption{ Average variance of the displacement error for the different filters in the case of uniform noise.}\label{F_unif}
\end{figure}

We consider the following estimators which are based on the aforementioned nominal state-space model:
\begin{itemize}
  \item KF denotes the standard Kalman filter;
  \item U-RKF denotes the proposed U-RKF with a  fixed tolerance $c = 0.5$;
  \item P-RKF+or denotes the prediction-resilient Kalman filter proposed in \cite{ROBUST_STATE_SPACE_LEVY_NIKOUKHAH_2013}; its  tolerance is  chosen in each realization through an \emph{oracle}; the latter has access to the true state trajectory and chooses the tolerance minimizing the mean squared filtering error.
  \item W-RKF+or denotes the Wasserstein distributionally robust filter  proposed in \cite{abadeh2018wasserstein} where the tolerance\footnote{In  \cite{abadeh2018wasserstein} the tolerance parameter is called Wasserstein radius.} is chosen in each realization by an oracle whose definition is the same as before;
  \item G-RKF+or denotes the globalized  robust filter proposed in \cite{10654520} whose parameter called ``targeted level of error''    is chosen in each realization by an oracle whose definition is the same as before;
  \item S-RKF denotes the sliding window variational outlier-robust filter proposed in \cite{9748015} with the parameter settings suggested as in Table 4 of \cite{9748015}.
\end{itemize}

\begin{table*}[t]
\centering
\caption{Mean and standard deviation (Std) of the running time of U-RKF and W-RKF across the four different scenarios.}
\small                
\setlength{\tabcolsep}{4pt} 
\renewcommand{\arraystretch}{1.1} 
\begin{tabular}{|c|c|c|c|c|c|}
\hline
                       &                & Sensor drift & Uniform noise & Outlier-cont. noise & Nonlinearity \\ \hline
\multirow{2}{*}{U-RKF} & Mean            & 0.0171s & 0.0172s & 0.0219s & 0.0171s \\ \cline{2-6}
                       & Std             & 0.0021s & 0.0004s & 0.0034s & 0.0006s \\ \hline
\multirow{2}{*}{W-RKF} & Mean            & 4.4711s & 4.4082s & 4.5328s & 4.4961s \\ \cline{2-6}
                       & Std             & 0.0297s & 0.0209s & 0.5910s & 0.1736s \\ \hline
\end{tabular}
\label{runningt}
\end{table*}

\begin{figure}[t]
  \centering
  \includegraphics[width=3.2in]{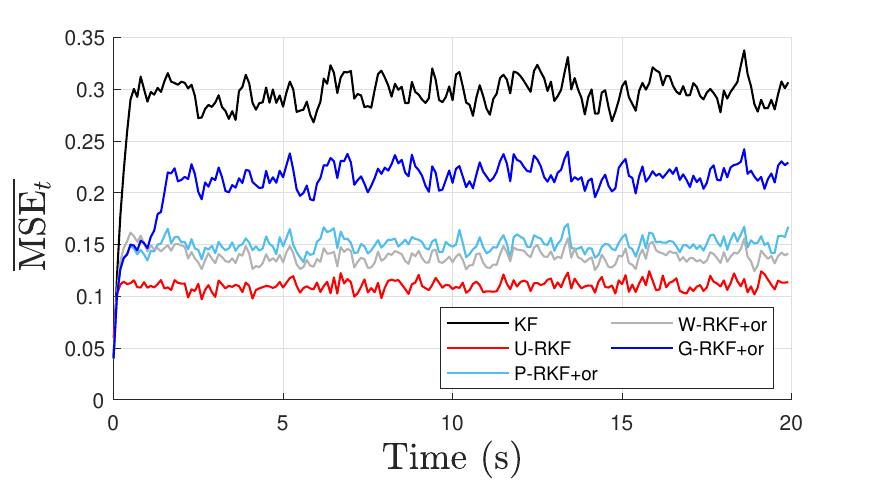}
  \caption{ Average variance of the displacement error for the different filters in the  presence of nonlinearities.}\label{F_out}
\end{figure}
\begin{figure}[t]
  \centering
  \includegraphics[width=3.2in]{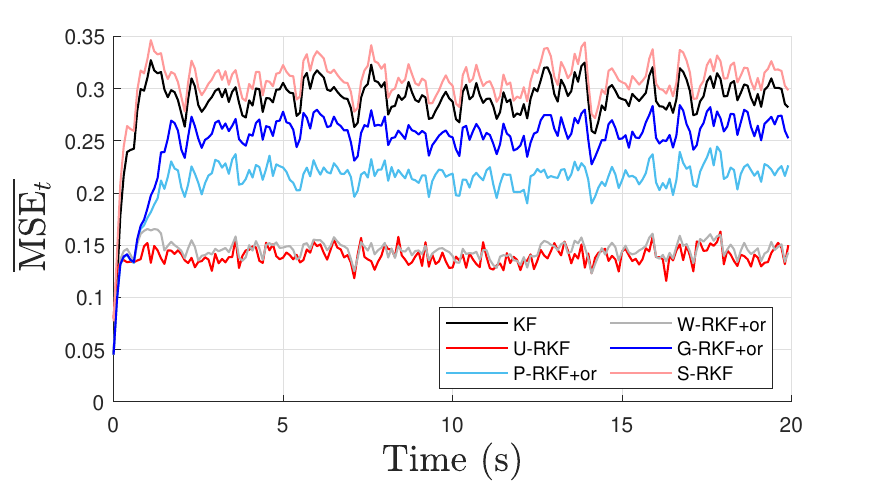}
  \caption{ Average variance of the displacement error for the different filters in the case of outlier-contaminated noise.}\label{F_dead}
\end{figure}

The oracles used in P-RKF+or,  W-RKF+or  and G-RKF+or search the optimal parameter (e.g. tolerance or targeted level of error which is constant over the time horizon) over a discretized interval of ten points whose extremal points are chosen in such a way the mean squared filtering error is properly captured without being overly broad, ensuring accuracy.

Then, we evaluate their performance through the average variance of the displacement error:
$\overline{\text{MSE}_t} = \frac{1}{M} \sum_{k=1}^{M}  \|\hat p^k_{t|t} -  p^k_t\|^2 $
where $ p^k_{t} $ and $\hat p^k_{t|t} $ denote the true displacement and  the estimated one corresponding to the $k$-th trajectory.  Fig. \ref{F_drif}–\ref{F_out} show the average variance of the displacement error for the estimators  across different scenarios.  It is worth stressing that the  ambiguity set for the proposed U-RKF is fixed a priori, whereas for P-RKF, W-RKF and G-RKF, it is provided by the oracle. Even so, as we can see, U-RKF achieves the best performance across  all three scenarios. The empirical evidence indicates that even when   a mild uncertainty is present in the state process, suggesting that one might naturally consider ambiguity sets like those used in P-RKF+or and G-RKF+or, it is ultimately more effective to adopt U-RKF, as it better balances robustness and performance.  In plain words, P-RKF+or and G-RKF+or are not able to split suitably the mismodeling budget between the state process and measurement models in this situation. As a consequence, these estimators are overly conservative in respect to the state process model and overly optimistic in respect to the observation model. Furthermore, W-RKF+or achieves a similar estimation performance to  U-RKF  because  the uncertainty  in both is framed in terms of the conditional density of $w_t$ given $Y_{t-1}$ as in  (\ref{minimax2psi}).  It is worth stressing that we did not report the performance of S-RKF in these cases because the comparison would not be fair, as S-RKF is specifically designed to handle the presence of outliers. Indeed, S-RKF performed worse than the robust estimators in all these cases. The performance of the estimators in the presence of outlier-contaminated noise is depicted in Fig. \ref{F_dead}. As shown, U-RKF outperforms P-RKF+or and G-RKF+or, and it performs similarly to W-RKF+or. Moreover, its performance is also better than that of S-RKF. This is because S-RKF is more competitive in scenarios where outliers also affect the state process model (\ref{nomi_linea_mod_proc}).

\begin{figure}[t]
  \centering
  \includegraphics[width=3.2in]{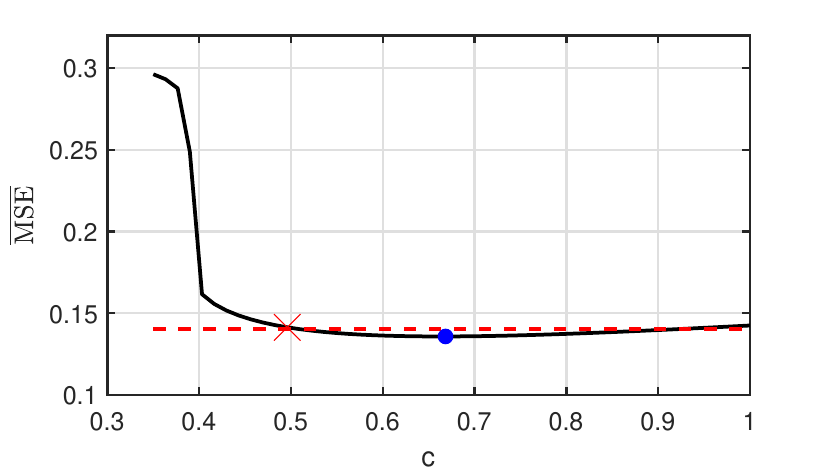}
  \caption{ $\overline{\text{MSE}}$ for U-RKF with $c \in  \mathcal{C}$ (solid black line), W-RKF+or (red dashed line). The red cross denotes the $\overline{\text{MSE}}$   for $c=0.5$ and the blue point the smallest value of $\overline{\text{MSE}}$.} \label{c_selec}
\end{figure}

To further assess how the choice of the tolerance $c$ influences the performance of the proposed robust estimator, we conducted the following experiment. Specifically, we consider a range of values
for $c$, that is $\mathcal{C} = [0.35,\, 1]$,
and assess the performance of U-RKF, evaluating the average variance of the displacement error on the entire time horizon 
$\overline{\mathrm{MSE}} = \frac{1}{N}\sum_{t=1}^{N} \overline{\mathrm{MSE}}_{t}.$
Figure~\ref{c_selec} reports the $\overline{\mathrm{MSE}}$ of U-RKF as a function of $c \in \mathcal{C}$ in the case of outlier-contaminated noise (black solid line). The red cross indicates the $\overline{\mathrm{MSE}}$ for U-RKF with $c = 0.5$ (i.e. the value considered in the previous analysis), while the blue point corresponds to the smallest value for the $\overline{\mathrm{MSE}}$. As shown, the performance of U-RKF does not significantly change for $c$ belonging to the interval $[0.5,0.8]$. Finally, similar results were observed for the other three cases.

Finally, we compare the computational time of U-RKF and W-RKF, i.e. the estimators showing the best performance in respect to the others. Here, all simulations have been implemented in MATLAB and  executed on a Mechanical Revolution notebook equipped with an $\text{AMD R9-7845HX CPU}$ and a $\text{GeForce  RTX 4070 GPU}$.  To ensure a fair comparison, we also apply a fixed tolerance for W-RKF, set equal to the mean of the tolerances provided by the oracle over the $M = 1000$ realizations.  In plain words, we consider a fair situation in which both algorithms compute the state estimate by considering only one tolerance. As shown in Table \ref{runningt}, the average running time per trial for U-RKF is approximately 0.02 seconds, while for W-RKF, it is around 4.5 seconds.
Moreover, the standard deviations (Std) of the running times for both U-RKF and W-RKF are significantly smaller than their respective means, indicating high stability and a low coefficient of variation. This demonstrates that the proposed U-RKF is significantly more computationally efficient than W-RKF.
The main reason is that W-RKF computes, at each time step, the filtering gain   using a gradient-like method, which is computationally expensive. In contrast, in our algorithm it is only required to compute a scalar quantity, i.e. $\theta_t$, which can be done using a bisection method, as done in \cite{RS_MPC_IET}.

\section{Update risk sensitive filter}\label{sec_rsf}
The minimax problem (\ref{minimax_psi}) can be relaxed. More precisely, the constraint on the maximizer can be replaced by  a penalty term in the objective function as follows:
\begin{equation}
 \label{minimax_rsk} \arg \underset{g_t \in \mathcal{G}_{t}} {\mathrm{min}}\max_{\tilde{\psi}_{t} \in \bar{\mathcal{B}}_{t}} J_t(\tilde {\psi}_t,g_t) - \theta \Ds(\tilde \psi_t, \psi_t)
\end{equation}
where $ \bar{\mathcal{B}}_{t}$ is the set of conditional densities of $y_t$ given $x_t$ and $\theta>0 $ is the risk sensitivity parameter which is a priori fixed.
Thus, there is no need to determine $\theta$ at every time step, thereby reducing the computational complexity.

\begin{thm}\label{th3}
Consider the estimation problem corresponding to the state space model (\ref{nomi_linea_mod_proc})-(\ref{nomi_linea_mod_obv}) and whose update estimate is obtained though (\ref{minimax_rsk}). Assume that the initial state is Gaussian distributed as in (\ref{pxY0}). Then, the  estimate of $x_t$ given $Y_t$ is  obtained by (\ref{hatxkf})-(\ref{Lt}) and  the corresponding error covariance matrix is
\begin{equation}\nn V_{t|t} = (P_{t|t}^{-1} - \theta I) ^{-1}.\end{equation}
where $P_{t|t}$ is  defined as in (\ref{ptt}), and $\theta$ must be such that $0<\theta<\sigma_{max}(P_{t|t})$. Moreover, condition (\ref{pxY2}) holds and the predictor takes the form (\ref{predxkf})-(\ref{predP}).
\end{thm}

\begin{pf}
Following the same reasonings in Lemma 1 and Theorem \ref{th1}, it is not difficult to prove that the least favorable density takes the form in (\ref{psi_0}) where $\lambda^{-1}_t$ is replaced by $\theta$. Then, it is possible to show that  the minimax problem  (\ref{minimax_rsk}) is  equivalent to
\begin{equation*}\label{minimax_r2}
 \hat x_{t|t} =\arg \underset{g_t \in \mathcal{G}_{t}}{\mathrm{min}}\max_{\tilde{p}_{t} \in { \check{\mathcal{B}}}_{t}}  \tilde  J_t (\tilde {p}_t,g_t) - \theta \Ds(\tilde p_t, p_t),
\end{equation*}
where  $\tilde J_t$ has been defined in (\ref{objt}), $\check{\mathcal{B}}_{t} $ is the set of  conditional densities of $ w_t := \begin{bmatrix} x_t^{\top} & y_t^{\top} \end{bmatrix}^{\top} $ given  $Y_{t-1}$.
Then, it is possible to show that  ${p}_{t}(w_t|Y_{t-1})$ and its maximizer, i.e. $\tilde{p}^0_{t}(w_t|Y_{t-1}),$ are both Gaussian. Thus,   we can deduce that $\tilde p_t(x_{t+1}|Y_{t})$ is Gaussian and the  corresponding  estimator.\qed
\end{pf}
The resulting estimator, which will be called update risk sensitive filter (U-RSF), is outlined in Algorithm  \ref{RSF}. Note that, U-RSF applies the distortion on $P_{t|t}$, which  is different  from   the  classic risk  sensitive filter \cite{H_INF_HASSIBI_SAYED_KAILATH_1999}, hereafter called prediction risk sensitive filter (P-RSF). The latter applies the distortion on $P_t$ because it assumes uncertainty is present in both (\ref{nomi_linea_mod_proc}) and (\ref{nomi_linea_mod_obv}).  
It remains  to characterize  the maximizer of (\ref{minimax_rsk}), i.e. the least favorable model corresponding  to U-RSF.
Following similar reasonings as the ones in the proof of Theorem \ref{th2}, it is possible to show that the least favorable model over the time interval    $\{ 0\ldots  N\}$ is the same to the one in (\ref{LF_x2}), expect that  \(\theta_t\) in (\ref{W}) is replaced  by  \(\theta\), i.e.
$W_{t+1}:=   \theta  I_n + \Omega^{-1}_{t+1}.$

\begin{algorithm}[t]
    \caption{ U-RSF at time step $t$}\label{RSF}
    \hspace*{\algorithmicindent} \textbf{Input} $\hat{x}_{t}$, $P_{t}$, $\theta$, $y_t$,  $ A$, $ Q$, $C$, $ R$ \\
    \hspace*{\algorithmicindent} \textbf{Output} $\hat{x}_{t|t}$, $\hat{x}_{t+1}$
    \begin{algorithmic}[1]
    \State $L_t = P_{t}  C^{\top}(CP_{t}C^{\top} + R)^{-1}$
    \State $\hat{x}_{t|t} = \hat{x}_{t}  + L_t(y_t -  C \hat{x}_{t})$
    \State $P_{t|t} = P_{t}- P_{t} C^{\top}( C P_{t} C^{\top} + R)^{-1}CP_{t}$
    \State $V_{t|t} = (P_{t|t}^{-1} - \theta I) ^{-1}$ \label{alg3_5}
    \State $\hat{x}_{t+1} =  A \hat{x}_{t|t}  $
    \State $P_{t+1} =  AV_{t|t}A^{\top}  + Q$
    \end{algorithmic}
\end{algorithm}

\subsection{Filter convergence}
We consider the situation  in which the pair $(A, C)$ is observable. The pair  $(A, Q)$ is reachable because $Q>0$. We show that it is possible to characterize an upper bound for the risk sensitivity parameter,  say $\theta_{MAX}$, which guarantees that the gain $L_t$ is well defined and converges to a constant value as $t\rightarrow \infty$ for any $\theta \in (0, \theta_{MAX}]$. By Algorithm \ref{RSF},  the covariance matrix $P_t$ obeys the recursion:
\al{\label{rec2} P_{t+1}=r^{RS}_\theta(P_t):=A(P_t^{-1}+C^\top R^{-1}C-\theta I)^{-1}A^\top+Q.}
Under the reachability and observability assumptions, there exists $\phi_k>0$, defined in the same way of the one in Section \ref{sec_4}, such that if $\theta \leq \phi_k$ then, by Proposition 3.1 in \cite{LEVY_ZORZI_RISK_CONTRACTION}, the sequence generated by (\ref{rec2}) converges and the corresponding algebraic equation admits a unique solution $P> 0$. However, unlike the recursion (\ref{rec}), it is not guaranteed that $P^{-1}_{t|t} - \theta I$ is positive definite for any $t$. Next, we identify the conditions on $P_0$ which guarantee that  $V_{t|t}>0$ for any $t\geq 0$. It is not difficult to see that the recursion (\ref{rec2}) can be written introducing an arbitrary observer gain matrix $G\in\mathbb R^{n\times m}$
\al{\label{recG}r^{RS}_\theta(P)=(&A-\alpha GC) (P^{-1}-\Psi_{\theta,\alpha})^{-1}(A-\alpha GC)^\top \nn\\
&-X_{\theta,\alpha,P}\Phi_{\theta,\alpha,P}^{-1}X_{\theta,\alpha,P}^\top+GRG^\top+Q}
where
\al{\Psi_{\theta,\alpha}&=(1-\alpha^2)C^\top R^{-1} C-\theta I\nn\\
X_{\theta,\alpha,P}&=\alpha(A-\alpha GC)(P^{-1}-\Psi_{\theta,\alpha})^{-1}C^\top-G \nn\\
\Phi_{\theta,\alpha,P}&=\alpha^2C(P^{-1}-\Psi_{\theta,\alpha})^{-1}C^\top+R\nn}
and $0<\alpha\leq 1$. Then, we consider the Lyapunov equation
\al{\label{lyap}\Sigma_{\rho,\alpha}=\rho^2(A-\alpha GC) \Sigma_{\rho,\alpha}(A-\alpha GC)^\top +GRG^\top+Q .}
Since $(A,C)$ is observable, we can choose $\alpha$ and $G$ such that $A-\alpha GC$ is a stable matrix. Let $r<1$ be the maximum among the  modules of the eigenvalues of $A-\alpha GC$. Then, for $1<\rho<r^{-1}$ matrix $\rho(A-\alpha GC)$ is a stable and the Lyapunov equation admits a unique solution. It is not difficult to see that the latter is also positive definite because $(A,Q)$ is reachable. The next result shows, for $\alpha$, $G$ and $\rho$ chosen as above, the conditions on $P_0$ which guarantees that $P^{-1}_{t|t} - \theta I$,  or equivalently $V_{t|t}$, is positive definite.
\begin{prop} \label{propRS}Let
\al{\beta_{\rho,\alpha}=\sigma_{min}\left(\frac{\rho^2-1}{\rho^2} \Sigma_{\rho,\alpha}^{-1} +(1-\alpha^2) C^\top R^{-1}C \right).\nn}
If $P_0$ for the iteration (\ref{rec2}) satisfies $0< P_0\leq \Sigma_{\rho,\alpha}$  and $0\leq \theta\leq \beta_{\rho,\alpha}$, then $0<P_t\leq \Sigma_{\rho,\alpha}$, and $V_{t|t}>0$ for any $t\geq 0$.
\end{prop}
\begin{pf}
First, we show that $P_t\leq \Sigma_{\rho,\alpha}$ implies that $V_{t|t}>0$. Condition $\theta \leq \beta_{\rho,\alpha}$ is equivalent to
\al{\frac{\rho^2-1}{\rho^2} \Sigma_{\rho,\alpha}^{-1} +(1-\alpha^2) C^\top R^{-1}C -\theta I\geq 0.\nn}
Since $\rho>1$ and $0<\alpha\leq 1$, it follows that
$ \Sigma_{\rho,\alpha}^{-1} + C^\top R^{-1}C -\theta I> 0.$
Since $P_t\leq \Sigma_{\rho,\alpha}$, it follows that
$ P_t^{-1} + C^\top R^{-1}C -\theta I> 0$
and thus
$V_{t|t}=(P_t^{-1} + C^\top R^{-1}C -\theta I)^{-1}> 0.$ Next, we prove that $r^{RS}(\Sigma_{\rho,\alpha})\leq \Sigma_{\rho,\alpha}$. Subtracting $r^{RS}_\theta(\Sigma_{\rho,\alpha})$ in (\ref{recG}) from (\ref{lyap}) we have
\al{&\Sigma_{\rho,\alpha}- r^{RS}_\theta(\Sigma_{\rho,\alpha}) =X_{\theta,\alpha,P} \Phi_{\theta,\alpha,P}^{-1}X_{\theta,\alpha,P}^\top  +(A-\alpha GC) \nn\\ & ~~~ \times (\rho^2 \Sigma_{\rho,\alpha}-(\Sigma_{\rho,\alpha}^{-1}-\Psi_{\theta, \alpha})^{-1} ) (A-\alpha GC)^\top \geq 0\nn}
because condition $\theta \leq \beta_{\rho,\alpha}$ implies the conditions
\al{\rho^2\Sigma_{\rho,\alpha}-(\Sigma_{\rho,\alpha}^{-1}-\Psi_{\theta, \alpha})^{-1}\geq 0, \quad \Phi_{\theta,\alpha,P}\geq 0.\nn}
Finally, we prove that if $P_0\leq \Sigma_{\rho,\alpha}$, then $P_t\leq \Sigma_{\rho,\alpha}$. By the monotonicity of the operator $r^{RS}_\theta$, see \cite[Lemma 5.1]{LEVY_ZORZI_RISK_CONTRACTION}, we have
$ P_1=r^{RS}_\theta(P_0)\leq r^{RS}_\theta(\Sigma_{\rho,\alpha})\leq \Sigma_{\rho,\alpha}.$
By induction, assume that $P_t\leq \Sigma_{\rho,\alpha}$ then
$ P_{t+1} = r^{RS}_\theta(P_t)\leq r^{RS}_\theta(\Sigma_{\rho,\alpha})\leq \Sigma_{\rho,\alpha}.$\qed
\end{pf}
We are ready to establish the convergence result.
\begin{thm}\label{teoRS2}
Let model (\ref{nomi_linea_mod_proc})-(\ref{nomi_linea_mod_obv}) be such that  $(A,C)$ is observable. Let $0<\alpha\leq 1$, $1<\rho<r^{-1}$,  $G\in \mathbb R^{n\times m}$, and $\phi_k$ with $k\geq n,$ chosen as before. We define
\al{\theta_{MAX}=\min\{\beta_{\rho,\alpha},\phi_k \}>0.\nn}
If $\theta\in(0,\theta_{MAX}]$, then  the sequence generated by the iteration (\ref{rec2}) converges to a unique solution $P>0$ for any $0<P_0 \leq \Sigma_{\rho,\alpha}$. Moreover, the limit $L$ of the filtering gain $L_t$ as $t\rightarrow \infty$ has the property that $A(I-LC)$ is stable.
\end{thm}
\begin{pf}
The convergence of $L_t$ follows by the previous reasonings, i.e. the combination of Proposition 3.1 in \cite{LEVY_ZORZI_RISK_CONTRACTION} and Proposition \ref{propRS}. The stabilizing property of the limit of $L_t$ can be proved likewise the proof of Theorem \ref{teo_c}.\qed
\end{pf}

In view of the above result, $\theta_{MAX}$ depends on $\rho,\alpha$ and $G$. Thus, it is possible to compute the best upper bound for $\theta$ optimizing  $\theta_{MAX}$ with respect to $(\rho,\alpha,G)$. It is worth noting that (\ref{rec2}) is the same recursion for the usual risk sensitive filter, i.e. P-RSF. However, there is a fundamental difference: while in Algorithm \ref{RSF} we require that
$\theta$ satisfies
$P_t^{-1}+C^\top R^{-1}C-\theta I >0,$
 in P-RSF we require the stronger condition
$ P_t^{-1}-\theta I >0.$
Thus, in general the range of $\theta$ in U-RSF is larger than the one for P-RSF. Moreover, the upper bound obtained through $\theta_{MAX}$ for the update case is in general larger than the prediction case. Indeed, forcing $\alpha=1$ in Theorem \ref{teoRS2} we obtain $\theta_{MAX}$ for the prediction case, see \cite{LEVY_ZORZI_RISK_CONTRACTION}.

Regarding the least favorable model in steady state,  we mention that it is possible to establish a result similar to  Proposition \ref{teoLFConv} where the condition on $c$ sufficiently small is  replaced by $\theta$  sufficiently small. Then, using the same reasoning at the end of Section \ref{sec_4}, it is possible to conclude that the covariance matrix of state estimation error  of U-RSF in steady state  is bounded.


\subsection{Numerical examples}
Consider  the nominal state space model   {\mz with the same constant matrices as those in (\ref{model1})}
and the initial state is modeled   as a Gaussian random vector
with zero mean and covariance matrix   $\mz P_0 = 0.1 I$.
First, since such model is both reachable and observable,  fixing $k=10$, we found $\mz \phi_k = 0.8387$.
In view of Theorem \ref{teoRS2}, we know that $\theta_{MAX}$ depends on  $\rho,~\alpha,~G$.
We maximize $\theta_{MAX}$ with $\alpha \in  (0,1]$, {\mz the entries of $G$ lie in the interval $[-10, 10]$} and $\rho \in (1, (\sigma_{max}(A-\alpha GC))^{-1}).$ We have found that the maximum value of $\theta_{MAX}$ for U-RSF {\mz is 0.5125.} Moreover, the maximum value of $\theta_{MAX} $ for P-RSF, computed in the same way as for U-RSF, but with $\alpha$ fixed equal to 1,  {\mz is 0.1709}.  We conclude that the upper bounds  for $\theta_{MAX} $ of U-RSF and P-RSF are  different and the one of U-RSF is larger than the one of P-RSF.

\begin{figure}
  \centering
  \includegraphics[width=3.2in]{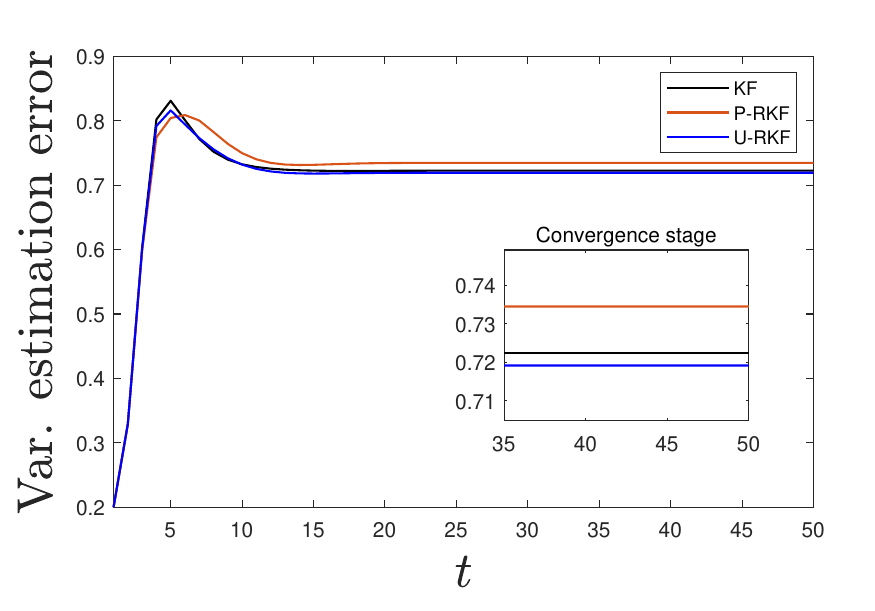}
  \caption{Variance of the estimation error when the standard KF (black line), P-RSF (red line) and U-RSF (blue line) are applied to the   least favorable model with $\mz \theta = 0.17$.}\label{F_6}
\end{figure}

Next, we compare the  performance of   U-RSF, P-RSF, and
KF applied to the  least favorable model which is the maximizer of (\ref{minimax_rsk}). We set $\mz \theta=0.17$ for  U-RSF, P-RSF and the least favorable model. Fig. \ref{F_6} shows that the proposed U-RSF performs better than P-RSF and the standard KF.  


Finally, we assessed the performance of U-RSF using the numerical experiments from Section~5.2. The risk-sensitivity parameter $\theta$ was chosen in each realization through an oracle, which has access to the true state trajectory and selects the value that minimizes the mean squared filtering error. We found that U-RSF performs better than KF but worse than the other estimators. This inferior performance can be explained as follows: since U-RSF is a risk-sensitive-like filter, it concentrates the uncertainty at a single time instant when the model is most susceptible to model uncertainty, which is not very realistic.

\section{Conclusion}\label{sec_con}
In this paper, we have studied a robust state estimation problem where the uncertainty is primarily in the  model for the observations. To address this, we have proposed a new robust estimation paradigm which is based on an ambiguity set that captures only the ``mismatch'' between the actual and nominal  observation models. The resulting robust estimator exhibits a structure similar to that of the Kalman filter  where the robustification takes place in the update stage. The latter is fundamentally different from  the  robust approaches  in \cite{ROBUST_STATE_SPACE_LEVY_NIKOUKHAH_2013,10654520} where the robustification takes place in the prediction stage. We have presented a numerical example based on a mass spring damper system, where sensor uncertainties constitute the dominant source of uncertainty. This example has shown that our estimator  outperforms the ones  in \cite{ROBUST_STATE_SPACE_LEVY_NIKOUKHAH_2013,10654520}. Furthermore, this example has shown that the robust estimator proposed in \cite{abadeh2018wasserstein} 	(equipped with an oracle) performs similarly to our approach.  Indeed, the robustification takes place in the update stage as in our estimator.  Thus, our analysis seems to suggest that the estimator in \cite{abadeh2018wasserstein} postulated uncertainty only in the observation model. Finally, the numerical results showed that our approach is preferable  since it is  significantly more computationally efficient than the one in \cite{abadeh2018wasserstein}.

\end{document}